\documentclass[preprint,12pt]{elsarticle}
\usepackage{natbib}

\usepackage{amsmath,amssymb,amsfonts}
\usepackage{amsthm}
\usepackage{mathrsfs}
\usepackage{mathtools}

\usepackage{graphicx}
\usepackage{multirow}
\usepackage{booktabs}
\usepackage{siunitx}
\usepackage{subcaption}
\usepackage{float}
\usepackage{threeparttable}
\usepackage{tabularx}

\usepackage{tikz}
\usepackage{tkz-graph}
\usetikzlibrary{arrows,shapes,patterns,positioning,calc,arrows.meta}
\usepackage{pgfplots}

\usepackage{algorithm}
\usepackage{algorithmicx}
\usepackage{algpseudocode}
\usepackage{listings}

\usepackage{xcolor}
\usepackage{textcomp}
\usepackage{pifont}

\usepackage[title]{appendix}
\usepackage{manyfoot}
\usepackage{hyperref}

\newtheorem{remark}{Remark}
\newtheorem{definition}{Definition}
\begin{document}

\begin{frontmatter}

\title{Beyond Linear Additive and Hill Functions: A General Logistic Reformulation 
       of Delay-Coupled Gene Regulatory Networks with Equilibrium Analysis, Hopf Bifurcation, and Lipschitz Stability}

\author[inst1]{Ismail Belgacem\corref{cor1}}
\ead{ismail.belgacem.81@gmail.com}
\address[inst1]{Mezaourou Ghazaouet, Tlemcen 13421, Algeria}
\cortext[cor1]{Corresponding author}

\begin{abstract}
Hill functions, the dominant tool for modeling gene regulatory networks, carry
fundamental limitations: when the cooperativity exponent takes a non-integer
value, as it routinely does when fitting experimental dose-response data,
derivatives diverge at the origin, complex-valued arithmetic silently corrupts
ODE trajectories, and the zero output at zero activation concentration traps
models irreversibly in off-states that no intrinsic dynamics can escape. This
paper uses  logistic functions based models that are globally
$C^\infty$, real-valued for all arguments, and strictly positive at zero
concentration. The substitution resolves all
three pathologies simultaneously while preserving the sigmoidal dynamics that
make Hill functions biologically appealing. Working with the delay-coupled
two-gene mutual-activation and self-repression network of Vinoth et al.\ as a
concrete model system, we analyze two logistic reformulations of the Hill-based
delay differential equations: a linear additive activation combined with
logistic self-repression, and a fully sigmoidal formulation in which both
activation and self-repression are governed by logistic functions. A closed-form
parameter matching relationship $\lambda = n/\theta$ is used by equating the
slopes of logistic and Hill functions at their half-maximal points, ensuring
local input--output equivalence. Closed-form parameter expressions for the weighted
logistic formulation, $\kappa_1 = 4g_A$, $\theta_B = g_A$, and
$\lambda_1 = \ln 3/g_A$, are derived by simultaneously matching basal expression
rates and local activation slopes to the original Hill-linear hybrid model. For the linear additive logistic formulation,
the unique biologically feasible equilibrium is computed at
$(167.96,\,164.22)$~nM; for the weighted logistic formulation it is computed at
$(144.46,\,139.99)$~nM, the reduction arising from saturation of the bounded
activation term. In the delay-free case ($\tau=0$), local asymptotic stability is established
for both formulations by proving that the Jacobian trace is strictly negative
for all positive parameters; the system then remains stable for delays
$\tau\in(0,\tau_c)$ and loses stability via Hopf bifurcation at the critical
delay $\tau_c$. Hopf bifurcation analysis, via the full two-equation
transcendental system solved numerically, yields critical self-repression
delays of $\tau_c \approx 1.19$~min for the linear additive formulation and
$\tau_c \approx 1.64$~min for the weighted logistic formulation, with
higher-order bifurcation sequences characterised numerically in each case;
direct numerical simulation of complex dynamics and extreme events in the
logistic reformulations is identified as an important direction for future work.
A quantitative comparison of the two formulations shows that replacing linear
additive activation with weighted logistic activation reduces the global
Lipschitz constant $L_F$ from $18.74$~min$^{-1}$ to $5.54$~min$^{-1}$, an
approximately 70\% reduction permitting substantially larger integration steps,
while the Jacobian Lipschitz constant $L_{DF}$ decreases from $0.24$ to
$0.084$~nM$^{-1}$\,min$^{-1}$. 
\end{abstract}


\end{frontmatter}

%%==========================================================================
\section{Introduction}
%%==========================================================================

Gene regulatory networks constitute the fundamental control architecture of
living cells, orchestrating the spatiotemporal patterns of gene expression that
underlie development, homeostasis, and adaptation to environmental change.
Mathematical modeling of these networks has become an indispensable tool for
understanding the dynamical behaviors that emerge from transcriptional
interactions --- oscillations, bistability, multistability, and
chaos~\cite{farcot2019chaos,belgacem2025glass} --- and for guiding the rational
design of synthetic gene circuits in biotechnology and medicine. At the heart of
nearly every such model sits the Hill function, a formulation introduced over a
century ago to describe cooperative ligand binding and subsequently adopted as
the near-universal representation of sigmoidal gene regulatory responses. Its
two forms, the activation function $h^+(x,\theta,n) = x^n/(x^n+\theta^n)$ and
the repression counterpart $h^-(x,\theta,n) = \theta^n/(x^n+\theta^n)$, are
intuitive, mechanistically grounded, and have been applied successfully to
models of the lambda phage lysis--lysogeny decision, the \textit{lac} and
\textit{gal} operons in \textit{Escherichia coli}, developmental patterning in
\textit{Drosophila}, and mammalian cell cycle control, among many others.
Complementary mathematical work has rigorously analyzed transcription-translation
systems and enzymatic reaction
networks~\cite{belgacem2018bmb,belgacem2014cdcRNA,%
belgacem2014med,belgacem2013cdc,belgacem2013cab,belgacem2013acta,belgacem2012ifac},
and has explored chaos, bifurcations, and hybrid dynamics in ring-circuit and
piecewise-linear gene network
architectures~\cite{farcot2019chaos,belgacem2025glass,belgacem2019ecc}.
The Hill coefficient $n$ quantifies cooperativity and $\theta$ the half-maximal
concentration, both carrying clear physical interpretations, and the sigmoidal
shape faithfully captures the switch-like regulatory transitions characteristic
of biological decision-making.

Yet this success conceals substantial hidden costs. The activation
Hill function vanishes identically at zero input, $h^+(0,\theta,n) = 0$. This
is mathematically convenient but biologically unrealistic. Experimental studies
consistently show minimal basal transcription (0.1--10\% of maximal) due to
stochastic promoter binding, incomplete repressor occupancy, and chromatin
dynamics. Models built on Hill functions either ignore this or introduce
ad~hoc additive offsets.

A second set of limitations concerns mathematical structure and numerical
behavior. When the Hill coefficient takes a non-integer value --- as it
routinely does when fitting experimental dose-response curves, yielding values
such as $n \approx 1.39$, $2.73$, or $3.52$ --- the function loses global
smoothness. For $n \in (k,k+1)$, derivatives of order greater than $k$ diverge
at the origin, restricting the function to $C^k$ smoothness and rendering
unavailable the higher-order analytical tools --- center manifold reduction,
normal form analysis, Hessian-based optimization --- that modern bifurcation
and parameter estimation workflows depend upon. The power-law form $x^n$ is
also computationally expensive for non-integer exponents, since it requires
evaluating the transcendental composition $e^{n\ln x}$, which accumulates
floating-point error near zero and, critically, becomes complex-valued whenever
any trajectory component transiently overshoots to a negative value, as every
adaptive ODE solver will produce from rounding errors alone. The result is
silent corruption: the solver continues integrating a complex-valued surrogate
system, producing smooth and visually plausible trajectories that are not
solutions to the true biological model. Beyond smoothness, the maximum slope of
$h^+$ at its inflection point equals $n/(4\theta)$, coupling threshold position
and transition steepness into a single expression. Threshold and steepness
cannot be adjusted independently, which obscures biological parameter
interpretation and complicates the rational tuning of synthetic circuits where
engineers need to set sensitivity and decision threshold as separate design
variables. Finally, the rational form $x^n/(x^n+\theta^n)$ has no closed-form
inverse for general $n$, blocking the feedback linearization and gradient-based
model predictive control
strategies~\cite{belgacem2020cdc,chambon2020automatica} that are increasingly
central to experimental implementation, while the zero production rate at $x=0$
creates fundamental controllability gaps from which Hill-based dynamics cannot
recover without external intervention.

The logistic function resolves all of these limitations simultaneously. Its
activation form $f^+(x,\theta,\lambda) = 1/(1+e^{-\lambda(x-\theta)})$ and
repression form $f^-(x,\theta,\lambda) = 1/(1+e^{\lambda(x-\theta)})$ are
globally $C^\infty$, real-valued for all arguments including negative ones,
and strictly positive at zero concentration for all finite $\lambda$ and
$\theta$. The elegant self-referential derivative $f' = \lambda f(1-f)$ reduces
Jacobian entries to products of function values, eliminating fractional
exponents from stability analysis entirely. The closed-form logit inverse
$f^{-1}(y) = \theta + \lambda^{-1}\ln(y/(1-y))$ enables exact feedback
linearization. The parameters $\theta$ and $\lambda$ are fully decoupled: the
threshold can be repositioned without altering the transition slope, and vice
versa, and both map directly to biologically measurable quantities --- $\theta$
to dissociation constants or half-maximal effective concentrations, and
$\lambda$ to effective cooperativity. Most importantly, the non-zero output at
$x = 0$ means that basal expression is built into the function's shape rather
than appended to it, with its magnitude controlled by the product $\lambda\theta$
through parameters that already define the regulatory response.

For more details, a preprint available on
arXiv~\cite{ismail2025logistic} that established the rigorous
mathematical foundations of the logistic alternatives to Hill functions:
global $C^\infty$ smoothness, real-valuedness for all arguments,
self-referential derivative structure
$\partial f^\pm/\partial x = \pm\lambda f^\pm(1-f^\pm)$, the closed-form
inverse via the logit transformation, and the uniform Lipschitz bound
$|\partial f^\pm/\partial x|\le\lambda/4$. The preprint also
proves the foundational existence-uniqueness-boundedness theorem for
the general product-of-logistics gene-regulatory-network framework.

The present manuscript builds on those foundations above but is independent
in its main results: it is concerned exclusively with the delay-coupled
two-gene network and develops, for the first time, the closed-form
parameter mappings between Hill-linear hybrid and pure logistic
formulations, the equilibrium and Hopf-bifurcation analysis of the
two reformulations, and the quantitative comparison of their
Lipschitz constants. 

As a concrete model system we take the delay-coupled two-gene
mutual-activation and self-repression network of Vinoth et
al.~\cite{vinoth2025extreme}, which exhibits periodic oscillations,
deterministic chaos, and extreme events under Hill-based formulation. We analyze
two logistic reformulations of this system: a linear additive formulation in
which the activation term $(g_A + g_{AB}B)$ is retained while self-repression
is replaced by a decreasing logistic, and a fully sigmoidal formulation in which
both activation and self-repression are logistic functions. For both
reformulations we use the parameter matching relationship $\lambda = n/\theta$
that equates logistic and Hill slopes at their respective half-maximal points,
compute unique biologically feasible equilibria, and prove local asymptotic
stability in the delay-free case by showing that the Jacobian trace is
strictly negative for all positive parameters. Hopf bifurcation analysis identifies critical
self-repression delays of $\tau_c \approx 1.19$~min for the linear additive
formulation and $\tau_c \approx 1.64$~min for the weighted logistic formulation,
and characterizes higher-order bifurcation sequences in each case; direct
numerical simulation of complex dynamics and extreme events in the logistic
reformulations is identified as important future work.

We quantify the reduction in Lipschitz constants --- approximately 70\% for
the vector-field constant $L_F$ and 65\% for the Jacobian constant $L_{DF}$
--- that follows from replacing unbounded linear activation with bounded
logistic forms, with direct consequences for integration step size and
numerical efficiency. We develop systematic parameter estimation methodologies
that convert hybrid Hill-linear models into pure logistic formulations by
simultaneously matching basal rates and activation slopes, deriving closed-form
expressions for all logistic parameters. We compare the linear additive and
weighted logistic formulations quantitatively. Together, these results establish
logistic functions as a principled, computationally robust, and biologically
faithful foundation for the next generation of gene regulatory network models.

%%==========================================================================
\section{The Logistic Function Framework}
\label{sec:logistic_framework}
%%==========================================================================

For the reader's convenience we recall the two logistic primitives
that underpin all subsequent analysis. The proofs of the analytical
properties stated in this section --- $C^\infty$ smoothness,
self-referential derivatives, the closed-form inverse, the closed-form
antiderivative, and the uniform Lipschitz bound --- are given in the
companion preprint~\cite{ismail2025logistic} and are not repeated
here. Only the identities used in the equilibrium and Hopf-bifurcation
analyses of Sections~\ref{sec:equilibrium}--\ref{sec:lipschitz} are
restated below.

\begin{definition}[Increasing logistic --- activation]
\label{def:fplus}
\begin{equation}
  f^+(x,\theta,\lambda) \;=\; \frac{1}{1+e^{-\lambda(x-\theta)}}.
  \label{eq:fplus}
\end{equation}
Models gene activation: off below threshold $\theta$, on above it,
with sharpness $\lambda>0$ and inflection point $f^+(\theta,\theta,\lambda)=\tfrac{1}{2}$.
\end{definition}

\begin{definition}[Decreasing logistic --- repression / NOT]
\label{def:fminus}
\begin{equation}
  f^-(x,\theta,\lambda) \;=\; \frac{1}{1+e^{\lambda(x-\theta)}}
  \;=\; 1-f^+(x,\theta,\lambda).
  \label{eq:fminus}
\end{equation}
Implements Boolean NOT: high regulator concentration suppresses expression.
\end{definition}

\begin{figure}[htbp]
  \centering
  \includegraphics[width=\linewidth]{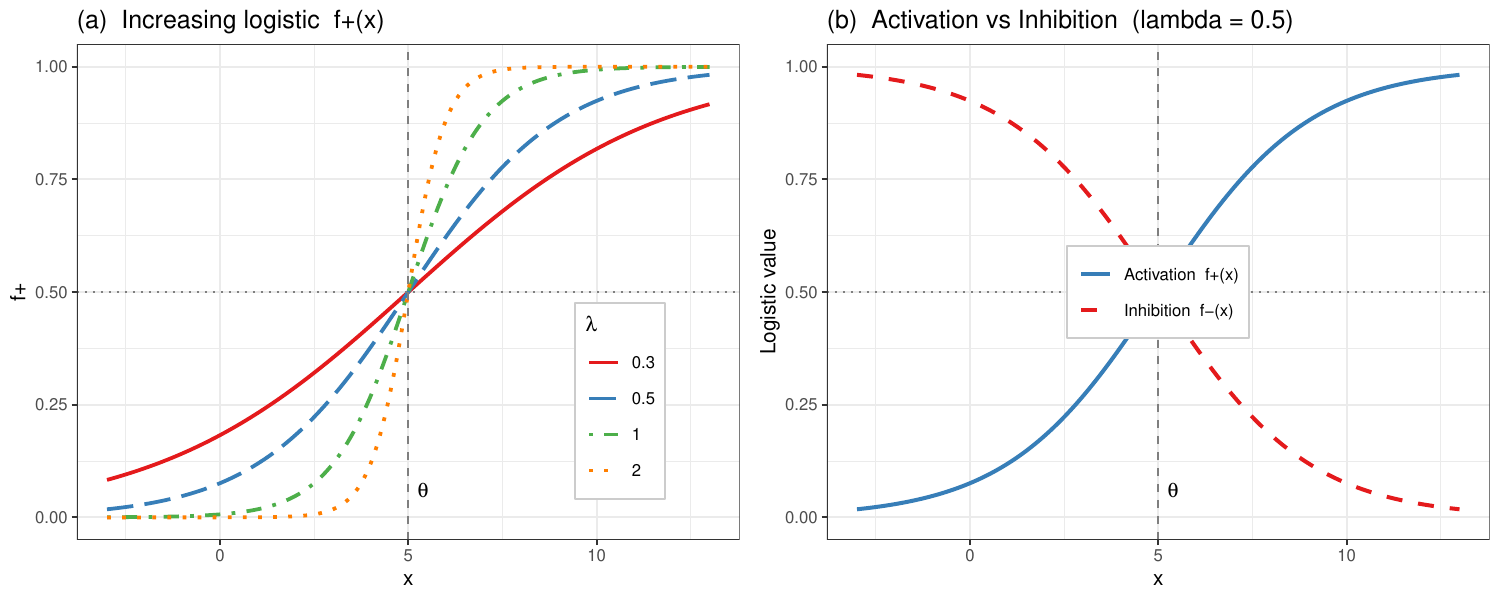}
  \caption{Logistic functions used in the Boolean logistic framework.
    \textbf{(a)} Increasing logistic $f^+(x,\theta,\lambda)$ for
    $\theta=5$ and $\lambda\in\{0.3,\,0.5,\,1.0,\,2.0\}$.
    \textbf{(b)} Comparison of $f^+$ and $f^-$ for $\theta=5$,
    $\lambda=0.5$; the dashed line marks $x=\theta$ where both
    attain half-maximum and maximum slope $\lambda/4$
    (Eq.~\eqref{eq:maxslope}).}
  \label{fig:logistic}
\end{figure}

Both functions are $C^\infty(\mathbb{R})$, strictly monotonic, and
bounded in $(0,1)$, as illustrated in Fig.~\ref{fig:logistic}.
Their self-referential derivative identities,
\begin{equation}
  \frac{df^+}{dx} = \lambda\,f^+(1-f^+) > 0, \qquad
  \frac{df^-}{dx} = -\lambda\,f^-(1-f^-) < 0,
  \label{eq:selfref}
\end{equation}
make Jacobian computation algebraic and enable all closed-form
linearizations developed throughout this paper.
The maximum slope of each function, attained at the inflection point
$x = \theta$, equals
\begin{equation}
  \max_x\left|\frac{df^\pm}{dx}\right| = \frac{\lambda}{4}.
  \label{eq:maxslope}
\end{equation}
Global $C^\infty$ smoothness contrasts sharply with the $C^k$ restriction
of non-integer Hill functions (for $n\in(k,k+1)$ the derivatives of
order greater than $k$ diverge at the origin); the full analytical
discussion of these Hill-function pathologies and their consequences for
numerical integration, and observer design appears
in~\cite{ismail2025logistic} and is not reproduced here.

%%==========================================================================
\section{Two-Gene Network: Model and Logistic Reformulations}
\label{sec:two_gene_model}
%%==========================================================================

\begin{figure}[h]
\centering
\begin{tikzpicture}[font=\small, >=stealth, node distance=3cm]
  \node[circle, draw, fill=yellow!20, minimum size=1.2cm, thick] (A) {$A$};
  \node[circle, draw, fill=yellow!20, minimum size=1.2cm, thick, right of=A] (B) {$B$};
  \draw[-|, thick, red] (A) .. controls +(140:1.5) and +(220:1.5) .. (A)
    node[midway, left, xshift=-0.1cm, font=\scriptsize] {self-repression};
  \draw[-|, thick, red] (B) .. controls +(40:1.5) and +(320:1.5) .. (B)
    node[midway, right, xshift=0.1cm, font=\scriptsize] {self-repression};
  \draw[->, thick, blue] (A.25) -- (B.155) node[midway, above=0.15cm, font=\scriptsize] {activation};
  \draw[->, thick, blue] (B.205) -- (A.335) node[midway, below=0.15cm, font=\scriptsize] {activation};
  \node[above left=0.15cm of A, font=\scriptsize] {$\tau_1$};
  \node[above right=0.15cm of B, font=\scriptsize] {$\tau_2$};
  \node[above=0.2cm of $(A)!0.5!(B)$, font=\scriptsize] {$\tau_{12}$};
  \node[below=0.2cm of $(A)!0.5!(B)$, font=\scriptsize] {$\tau_{21}$};
\end{tikzpicture}
\caption{Regulatory topology of the two-gene network. Red lines with bars
  indicate cooperative self-repression; blue arrows represent mutual activation.
  Delays $\tau_1, \tau_2$ govern self-repression lags, while $\tau_{12},
  \tau_{21}$ represent cross-activation delays.}
\label{fig:two_gene_regulation}
\end{figure}
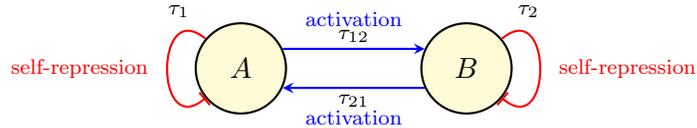

To illustrate the logistic framework concretely, we adopt the delay-coupled
two-gene mutual-activation and self-repression network studied by Vinoth et
al.~\cite{vinoth2025extreme}. The system exhibits complex behaviors including
periodic oscillations, deterministic chaos, and \emph{extreme events} ---
abrupt, large-amplitude spikes in protein concentration reminiscent of ``rogue
waves'' in nonlinear dynamics. The network comprises two proteins, $A$ and $B$,
with the following regulatory architecture:
\begin{itemize}
    \item \textbf{Mutual activation:} Protein $A$ activates gene $B$, and
          protein $B$ activates gene $A$.
    \item \textbf{Cooperative self-repression:} Each gene represses its own
          expression via multimeric binding (e.g., dimers or tetramers),
          traditionally captured by a Hill coefficient $n > 1$.
    \item \textbf{Time delays:} The system incorporates delays
          $\tau_1, \tau_2, \tau_{12}, \tau_{21}$ representing
          transcription/translation lags or signalling delays.
\end{itemize}
The regulatory architecture is depicted schematically in
Fig.~\ref{fig:two_gene_regulation}.

%----------------------------------------------------------------------
\subsection{Original Hill-Based Formulation}
\label{subsec:hill_form}
%----------------------------------------------------------------------

In the original study~\cite{vinoth2025extreme}, the model employs additive
linear activation combined with Hill-type repression. The delay differential
equations (DDEs) for protein concentrations $A(t)$ and $B(t)$ are:
\begin{align}
\frac{dA(t)}{dt} &= \left( g_A + g_{AB} B(t - \tau_{12}) \right)
  \frac{1}{1 + \left( \dfrac{A(t - \tau_1)}{A_0} \right)^n}
  - \gamma_A A(t), \label{eq:A_hill}\\[6pt]
\frac{dB(t)}{dt} &= \left( g_B + g_{BA} A(t - \tau_{21}) \right)
  \frac{1}{1 + \left( \dfrac{B(t - \tau_2)}{B_0} \right)^n}
  - \gamma_B B(t), \label{eq:B_hill}
\end{align}
where:
\begin{itemize}
    \item $g_A, g_B$ are basal production rates,
    \item $g_{AB}, g_{BA}$ quantify the strength of mutual activation,
    \item $\gamma_A, \gamma_B$ are degradation rates,
    \item $A_0, B_0$ are repression thresholds,
    \item $n$ is the Hill coefficient encoding cooperativity in self-repression.
\end{itemize}

All parameter values follow~\cite{vinoth2025extreme}: basal production
$g_A = g_B = 50$~nM/min, cross-activation strengths
$g_{AB} = g_{BA} = 3$~min$^{-1}$, degradation rates
$\gamma_A = 0.20$~min$^{-1}$ and $\gamma_B = 0.24$~min$^{-1}$, repression
thresholds $A_0 = B_0 = 100$~nM, and Hill coefficient $n = 4$ representing
tetrameric cooperative binding. 
%----------------------------------------------------------------------
\subsection{Linear Additive Logistic Reformulation}
\label{subsec:linear_reform}
%----------------------------------------------------------------------

The first reformulation retains the linear additive activation term
$(g_A + g_{AB}B)$ and replaces only the Hill self-repression term
$1/(1+(A/A_0)^n)$ by the decreasing logistic $f^-(A,A_0,\lambda_A)$.
The reformulated DDEs are:
\begin{align}
\dot{A}(t) &= \left( g_A + g_{AB} B(t - \tau_{12}) \right)
  f^-\!\big(A(t-\tau_1), A_0, \lambda_A\big) - \gamma_A A(t),
  \label{eq:A_logistic}\\[6pt]
\dot{B}(t) &= \left( g_B + g_{BA} A(t - \tau_{21}) \right)
  f^-\!\big(B(t-\tau_2), B_0, \lambda_B \big) - \gamma_B B(t).
  \label{eq:B_logistic}
\end{align}
Written out explicitly:
\begin{align}
\dot{A} &= \left( g_A + g_{AB} B(t - \tau_{12}) \right)
  \frac{1}{1+e^{\lambda_A (A(t-\tau_1)-A_0)}} - \gamma_A A(t),
  \label{eq:A_explicit}\\
\dot{B} &= \left( g_B + g_{BA} A(t - \tau_{21}) \right)
  \frac{1}{1+e^{\lambda_B (B(t-\tau_2)-B_0)}} - \gamma_B B(t).
  \label{eq:B_explicit}
\end{align}
Here $\lambda_A$ and $\lambda_B$ are steepness parameters derived via the
slope-matching criterion of Section~\ref{subsubsec:steepness_match}:
$\lambda_A = n/A_0$ and $\lambda_B = n/B_0$.

%----------------------------------------------------------------------
\subsection{Weighted Logistic Reformulation}
\label{subsec:weighted_reform}
%----------------------------------------------------------------------

The second reformulation replaces the Hill repression term \emph{and} the
unbounded linear activation term $(g_A + g_{AB}B)$ with bounded increasing
logistic functions. The maximal production rate $\kappa_1$ multiplies the
weighted-signal sigmoid:
\begin{align}
\dot{A}(t)
  &= \kappa_1\,
     f^+\!\bigl(g_{AB}B(t-\tau_{12}),\,\theta_B,\,\lambda_1\bigr)\,
     f^-\!\bigl(A(t-\tau_1),\,A_0,\,\lambda_3\bigr)
     - \gamma_A A(t),
  \label{eq:A_logistic_basal}\\[4pt]
\dot{B}(t)
  &= \kappa_2\,
     f^+\!\bigl(g_{BA}A(t-\tau_{21}),\,\theta_A,\,\lambda_2\bigr)\,
     f^-\!\bigl(B(t-\tau_2),\,B_0,\,\lambda_4\bigr)
     - \gamma_B B(t),
  \label{eq:B_logistic_basal}
\end{align}
or written explicitly:
\begin{align}
\dot{A}(t)
  &= \frac{\kappa_1}
          {1+e^{-\lambda_1\bigl(g_{AB}B(t-\tau_{12})-\theta_B\bigr)}}
     \cdot
     \frac{1}{1+e^{\,\lambda_3(A(t-\tau_1)-A_0)}}
     - \gamma_A A(t),
  \label{eq:A_explicit_basal}\\[4pt]
\dot{B}(t)
  &= \frac{\kappa_2}
          {1+e^{-\lambda_2\bigl(g_{BA}A(t-\tau_{21})-\theta_A\bigr)}}
     \cdot
     \frac{1}{1+e^{\,\lambda_4(B(t-\tau_2)-B_0)}}
     - \gamma_B B(t).
  \label{eq:B_explicit_basal}
\end{align}
The thresholds $\theta_B$ and $\theta_A$ carry units of nM\,min$^{-1}$
because they are compared against the weighted signals $g_{AB}B$ and
$g_{BA}A$, which have those units.

A direct calculation shows that, when the activator $B$ is set to zero
in~\eqref{eq:A_explicit_basal}, the production term reduces to
\begin{equation}
    \frac{\kappa_1}{1+e^{\lambda_1\theta_B}}\,
    f^-(A,A_0,\lambda_3),
    \label{eq:Abasal}
\end{equation}
which remains strictly positive for every finite choice of $\lambda_1$
and $\theta_B$. The reformulated system therefore inherits a built-in
floor on the production rate of gene $A$ — a property that is
quantitatively consistent with the large basal rate
$g_A = 50$~nM/min in the original Hill-based
model~\cite{vinoth2025extreme}, and that captures the constitutively
active behaviour of the promoter without any auxiliary additive offset.

\begin{remark}[Tunability of the basal regime via the product $\lambda_1\theta_B$]
\label{rem:basal_constitutive}
The single parameter combination $\lambda_1\theta_B$ controls how far
the basal expression~\eqref{eq:Abasal} sits below the saturation
ceiling $\kappa_1$. Large values of the product drive the prefactor
$1/(1+e^{\lambda_1\theta_B})$ towards zero, recovering tightly
repressed promoter behaviour; small values keep the prefactor near
$\tfrac{1}{2}$, reproducing constitutively active promoters such as
those in~\cite{vinoth2025extreme}. A practical consequence for
parameter identification is that a single basal-rate measurement at
$B=0$ provides one scalar constraint linking $\kappa_1$, $\lambda_1$
and $\theta_B$; combined with steady-state and transient measurements
this constraint resolves the otherwise-coupled identification of the
three parameters. The standard Hill--linear hybrid form, in which
$g_A$ enters as an additive constant outside the Hill factor, does
not admit such a single-observation anchor.
\end{remark}

%%==========================================================================
\section{Parameter Matching and Validation}
\label{sec:param_match_strategy}
%%==========================================================================

To ground our logistic-based modeling framework in biological reality,
we establish a systematic parameter estimation strategy that links
the logistic parameters to experimentally accessible quantities.
Our approach leverages time-series measurements of gene expression,
encompassing maximal production rates $\kappa_i$, degradation rates
$\gamma_i$, steepness parameters $\lambda_{ij}$, and regulatory
thresholds $\theta_{ij}$.  The smooth, bounded derivatives of the
logistic framework facilitate gradient-based optimization across
high-dimensional parameter landscapes.  Unlike Hill-based models,
where fractional exponents lead to discontinuous parameter
sensitivities, the exponential structure of logistic functions ensures
that small parameter perturbations yield proportionally small changes
in model predictions, enabling robust parameter identifiability even
from noisy experimental data.

When time-series data are available, all parameters can be estimated
by minimising the sum of squared residuals between observed gene
expression trajectories and model predictions:
\begin{equation}
\min_{\kappa_i,\,\gamma_i,\,\lambda_{ij},\,\theta_{ij}}
\sum_{t}
\left\|\mathbf{x}(t)
- \hat{\mathbf{x}}\!\left(t;\,\kappa_i,\gamma_i,\lambda_{ij},\theta_{ij}\right)
\right\|_2^2,
\label{eq:least_squares}
\end{equation}
where $\mathbf{x}(t)$ represents the experimental measurements at time
$t$ and $\hat{\mathbf{x}}(t;\,\cdot)$ denotes the model-predicted
expression levels.  The logistic function's global smoothness and
analytically tractable derivatives,
\begin{equation}
f'(s) = \frac{e^{-s}}{(1+e^{-s})^2}, \qquad
f''(s) = \frac{e^{-s}(e^{-s}-1)}{(1+e^{-s})^3},
\label{eq:logistic_derivatives}
\end{equation}
ensure that gradient-based algorithms such as Gauss--Newton or
Levenberg--Marquardt converge efficiently, in sharp contrast to
Hill-function models where fractional exponents and potential
singularities near the origin can create ill-conditioned landscapes.
In the absence of time-series data, the analytic matching strategy
developed below provides closed-form starting values for all logistic
parameters directly from the original Hill--linear hybrid model.

%----------------------------------------------------------------------
\subsection{Steepness Matching: $\lambda = n/\theta$}
\label{subsubsec:steepness_match}
%----------------------------------------------------------------------

To ensure local input--output equivalence between the Hill repression
term and its logistic replacement, we equate their slopes at the
half-maximal point $x = \theta = A_0$.

For the Hill repression function $h^-(x,\theta,n) = \theta^n/(x^n+\theta^n)$,
differentiation via the quotient rule gives
\begin{equation}
  \frac{dh^-}{dx} = \frac{-n\theta^n x^{n-1}}{(x^n+\theta^n)^2}.
\end{equation}
At $x=\theta$:
\begin{equation}
  \left.\frac{dh^-}{dx}\right|_{x=\theta}
  = \frac{-n\theta^{2n-1}}{4\theta^{2n}} = -\frac{n}{4\theta}.
\end{equation}
For the decreasing logistic $f^-(x,\theta,\lambda)$,
equations~\eqref{eq:selfref} and~\eqref{eq:maxslope} give the slope
at $x=\theta$:
\begin{equation}
  \left.\frac{df^-}{dx}\right|_{x=\theta} = -\frac{\lambda}{4}.
\end{equation}
Setting the two slopes equal:
\begin{equation}
  \frac{\lambda}{4} = \frac{n}{4\theta}
  \quad\Longrightarrow\quad
  \boxed{\lambda = \frac{n}{\theta}.}
  \label{eq:lambda_match}
\end{equation}
Applied to the self-repression terms of the two-gene network:
\begin{equation}
  \lambda_A = \frac{n}{A_0} = \frac{4}{100} = 0.04~\text{nM}^{-1},
  \qquad
  \lambda_B = \frac{n}{B_0} = \frac{4}{100} = 0.04~\text{nM}^{-1}.
  \label{eq:lambda_computed}
\end{equation}

%----------------------------------------------------------------------
\subsection{Slope and Basal Matching for the Weighted Logistic}
\label{subsubsec:slope_basal_match}
%----------------------------------------------------------------------

We derive the weighted logistic parameters by simultaneously matching
the \emph{basal expression rate} and the \emph{local activation slope}
of the original linear additive model.  Working with the weighted signal
$u = g_{AB}B$, the weighted logistic activation reads
$\kappa_1 f^+(u,\theta_B,\lambda_1)$.

We begin with the basal matching condition.
Unlike Hill function models, which require arbitrary offsets or ad~hoc
modifications to reproduce basal transcription, the logistic function
naturally encodes a non-zero baseline expression through its
mathematical structure.  Setting $u=0$ (i.e., $B=0$) and assuming
negligible self-repression, the logistic basal rate must equal the
linear model's basal rate $g_A$:
\begin{equation}
\frac{\kappa_1}{1 + e^{\lambda_1\theta_B}} = g_A.
\label{eq:basal_cond}
\end{equation}
Rearranging:
\begin{equation}
\lambda_1
= \frac{1}{\theta_B}\ln\!\left(\frac{\kappa_1}{g_A} - 1\right).
\label{eq:lambda1_from_basal}
\end{equation}
This relationship enables a bidirectional inference strategy: given
experimentally measured basal expression rates and estimated threshold
concentrations, $\lambda_1$ can be computed directly; conversely, if
cooperativity estimates are available from binding studies, threshold
positions can be inferred from basal measurements.

With the basal condition in place, we turn to slope matching via Taylor expansion.
For moderate steepness $\lambda_1$, the logistic activation admits a
first-order Taylor expansion around its inflection point $u = \theta_B$:
\begin{equation}
\kappa_1 f^+(g_{AB}B,\theta_B,\lambda_1)
\approx
\frac{\kappa_1}{2}
+ \frac{\kappa_1\lambda_1}{4}
  \bigl(g_{AB}B-\theta_B\bigr).
\label{eq:taylor_approx}
\end{equation}
Matching to the original term $g_A + g_{AB}B$ imposes two conditions:
\begin{enumerate}
\item \textbf{Slope condition}: the coefficient of $B$ must equal $g_{AB}$,
\begin{equation}
\frac{\kappa_1\lambda_1 g_{AB}}{4} = g_{AB}
\quad\Longrightarrow\quad
\kappa_1\lambda_1 = 4.
\label{eq:slope_cond}
\end{equation}
Note that $g_{AB}$ cancels; this condition is independent of the
activation strength.
\item \textbf{Intercept condition}: substituting $\kappa_1\lambda_1 = 4$
from~\eqref{eq:slope_cond} into the constant term,
\begin{equation}
\frac{\kappa_1}{2} - \theta_B = g_A
\quad\Longrightarrow\quad
\kappa_1 = 2(g_A + \theta_B).
\label{eq:intercept_cond}
\end{equation}
\end{enumerate}

The two matching conditions~\eqref{eq:slope_cond} and~\eqref{eq:intercept_cond}
leave the threshold $\theta_B$ free; we close the system with a biologically
motivated threshold selection.
We adopt the criterion that the sigmoid's inflection point coincides
with the weighted signal level at which cross-activation equals basal
expression:
\begin{equation}
g_{AB}B^* = g_A
\quad\Longrightarrow\quad
\theta_B = g_A.
\label{eq:threshold_criterion}
\end{equation}
This choice has three compelling properties.  First, it centers the
steepest sigmoid region at the biological transition between
basal-dominated and activator-driven expression (the inflection
displacement $g_{AB}B^* - \theta_B = 0$ vanishes exactly at this
transition).  Second, it guarantees $\theta_B/g_A = 1$,
validating the moderate-ratio assumption underlying the Taylor
approximation above.  Third, $\theta_B/g_{AB} = g_A/g_{AB}$ is
directly interpretable as an activation midpoint analogous to the
EC$_{50}$ in pharmacology, making it accessible to independent
dose-response measurements.

These three conditions together yield closed-form parameter expressions.
Substituting $\theta_B = g_A$ into~\eqref{eq:intercept_cond}:
\begin{equation}
\kappa_1 = 4g_A.
\label{eq:kappa1_final}
\end{equation}
Substituting $\kappa_1 = 4g_A$ and $\theta_B = g_A$
into~\eqref{eq:lambda1_from_basal}:
\begin{equation}
\lambda_1
= \frac{1}{g_A}\ln\!\left(\frac{4g_A}{g_A} - 1\right)
= \frac{\ln 3}{g_A}.
\label{eq:lambda1_final}
\end{equation}
By the complete symmetry of the network:
$\kappa_2 = 4g_B$, $\theta_A = g_B$, $\lambda_2 = \ln 3/g_B$.
For the self-repression terms, the steepness matching of
Section~\ref{subsubsec:steepness_match} gives
$\lambda_3 = n/A_0$ and $\lambda_4 = n/B_0$.
The complete parameter correspondence for the weighted logistic
formulation is:
\begin{align}
\kappa_1 &= 4g_A, & \kappa_2 &= 4g_B, \notag\\[4pt]
\theta_B  &= g_A, & \theta_A  &= g_B, \notag\\[4pt]
\lambda_1 &= \frac{\ln 3}{g_A}, & \lambda_2 &= \frac{\ln 3}{g_B}, \notag\\[4pt]
\lambda_3 &= \frac{n}{A_0},     & \lambda_4 &= \frac{n}{B_0}.
\label{eq:weighted_params}
\end{align}

These values, together with the derived logistic parameters, are summarised in Table~\ref{tab:parameters}.

\begin{table}[htbp!]
\centering
\caption{Parameter values for the two-gene regulatory network. Parameters
  $g_A$ through $n$ are drawn from Vinoth et al.~\cite{vinoth2025extreme};
  steepness parameters $\lambda_A$, $\lambda_B$ are derived via
  $\lambda = n/\theta$ (Section~\ref{subsubsec:steepness_match});
  weighted logistic parameters $\kappa_1$, $\kappa_2$, $\theta_B$,
  $\theta_A$, $\lambda_1$, $\lambda_2$ follow from
  equations~\eqref{eq:kappa1_final}--\eqref{eq:weighted_params}.}
\label{tab:parameters}
\begin{tabular}{llll}
\toprule
\textbf{Parameter} & \textbf{Symbol} & \textbf{Value} & \textbf{Units} \\
\midrule
Basal production (A)             & $g_A$       & 50    & nM\,min$^{-1}$ \\
Basal production (B)             & $g_B$       & 50    & nM\,min$^{-1}$ \\
Cross-activation (A$\to$B)       & $g_{AB}$    & 3     & min$^{-1}$ \\
Cross-activation (B$\to$A)       & $g_{BA}$    & 3     & min$^{-1}$ \\
Degradation rate (A)             & $\gamma_A$  & 0.20  & min$^{-1}$ \\
Degradation rate (B)             & $\gamma_B$  & 0.24  & min$^{-1}$ \\
Repression threshold (A)         & $A_0$       & 100   & nM \\
Repression threshold (B)         & $B_0$       & 100   & nM \\
Hill coefficient                 & $n$         & 4     & dimensionless \\
Steepness parameter (A)          & $\lambda_A$ & 0.04  & nM$^{-1}$ \\
Steepness parameter (B)          & $\lambda_B$ & 0.04  & nM$^{-1}$ \\
\midrule
Maximal production (A)           & $\kappa_1$  & 200   & nM\,min$^{-1}$ \\
Maximal production (B)           & $\kappa_2$  & 200   & nM\,min$^{-1}$ \\
Activation threshold (B signal)  & $\theta_B$  & 50    & nM\,min$^{-1}$ \\
Activation threshold (A signal)  & $\theta_A$  & 50    & nM\,min$^{-1}$ \\
Activation steepness (A eq.)     & $\lambda_1$ & $\approx 0.022$ & (nM\,min$^{-1}$)$^{-1}$ \\
Activation steepness (B eq.)     & $\lambda_2$ & $\approx 0.022$ & (nM\,min$^{-1}$)$^{-1}$ \\
Repression steepness (A eq.)     & $\lambda_3$ & 0.04  & nM$^{-1}$ \\
Repression steepness (B eq.)     & $\lambda_4$ & 0.04  & nM$^{-1}$ \\
\bottomrule
\end{tabular}
\end{table}

%----------------------------------------------------------------------
\subsection{Formal Equivalence: Weighted and Fixed-Weight Product-of-Logistics}
\label{subsubsec:formal_equiv}
%----------------------------------------------------------------------
Within the product-of-logistics framework, incorporating explicit
positive interaction weights $w_{ij} > 0$ is formally equivalent to
the fixed-weight formulation after a parameter rescaling.  For an
activator with weight $w_{ij} > 0$:
\begin{equation}
f^+(w_{ij}x_j,\,\theta_{ij},\,\lambda_{ij})
= \frac{1}{1+e^{-\lambda_{ij}(w_{ij}x_j - \theta_{ij})}}
= \frac{1}{1+e^{-\lambda'_{ij}(x_j - \theta'_{ij})}},
\label{eq:weighted_equiv}
\end{equation}
where $\lambda'_{ij} = \lambda_{ij}w_{ij} > 0$ and
$\theta'_{ij} = \theta_{ij}/w_{ij} > 0$.  The effective threshold
$\theta' = \theta/w > 0$ is always positive and biologically
interpretable as EC$_{50}$.  An identical rescaling holds for
decreasing logistic repression terms.  Consequently, the weighted and
fixed-weight formulations produce identical trajectories after
parameter estimation and differ only in how steepness and threshold
information is distributed across parameters.

\begin{remark}[Parameter units: weighted vs.\ fixed-weight formulations]
\label{rem:param_comparison}
The two logistic reformulations use the same parameter names but with
different physical units, because their logistic arguments have
different dimensions.

\begin{center}
\renewcommand{\arraystretch}{1.3}
\begin{tabular}{lllll}
\toprule
Formulation & Logistic argument & $\lambda_1$ value
  & $\lambda_1$ units & $\theta_B$ units \\
\midrule
Fixed-weight  & $B - \theta_B$
  & $\approx 0.066$ & $\mathrm{nM^{-1}}$ & nM \\
Weighted & $g_{AB}B - \theta_B$
  & $\approx 0.022$ & $(\mathrm{nM\,min^{-1}})^{-1}$ & $\mathrm{nM\,min^{-1}}$ \\
\bottomrule
\end{tabular}
\end{center}

The fixed-weight value $\lambda_1' \approx 0.066~\mathrm{nM^{-1}}$ is
obtained from the rescaling $\lambda_1' = \lambda_1\,g_{AB}
= (\ln 3/g_A)\,g_{AB} = (\ln 3/50)\times 3 \approx 0.066$, using
the main parameter set ($g_A = 50~\mathrm{nM\,min^{-1}}$,
$g_{AB} = 3~\mathrm{min^{-1}}$).

Both formulations produce nearly identical trajectories
(Fig.~\ref{fig:logistic_compare_vinoth_2}) because the dimensionless
products $\lambda_1(B-\theta_B)$ and $\lambda_1(g_{AB}B - \theta_B)$
take the same numerical values at biologically relevant concentrations
when parameters are chosen via the matching procedure above.  The
choice between formulations should be guided by which parameterisation
maps more directly to available experimental measurements.
\end{remark}

\begin{remark}[Approximate nature of the parameter matching]
\label{rem:approx_matching}
The slope condition~\eqref{eq:slope_cond} and the exact basal
condition~\eqref{eq:basal_cond} cannot be satisfied simultaneously.
With $\theta_B = g_A$ and $\kappa_1 = 4g_A$, the formula
$\lambda_1 = \ln 3/g_A$ enforces \emph{exact} basal matching.
A consistency check shows that the implied slope coefficient at the
sigmoid midpoint is
\begin{equation}
\frac{\kappa_1\lambda_1}{4}
= \frac{4g_A \cdot \tfrac{\ln 3}{g_A}}{4}
= \ln 3 \approx 1.099,
\end{equation}
a $\approx 10\%$ overestimate of the exact value of unity required
by~\eqref{eq:slope_cond}.  For moderate steepness regimes this
discrepancy is negligible, as confirmed by the close trajectory
agreement in Fig.~\ref{fig:logistic_compare_vinoth_2}.
\end{remark}

%----------------------------------------------------------------------
\subsection{Alternative Threshold Selection Strategies}
\label{subsubsec:alt_threshold}
%----------------------------------------------------------------------

While the biologically motivated threshold selection
$\theta_B = g_A$, $\theta_A = g_B$ (weighted formulation) provides
elegant closed-form parameter expressions with clear biological
interpretation, alternative strategies merit consideration depending
on data availability and modeling objectives.

One natural alternative is to set thresholds near the system's
steady-state concentrations, which can be determined experimentally
or computed numerically.  For the original system with linear
activation, the steady-state equations are
\begin{align}
\gamma_A A^* &= \bigl(g_A + g_{AB}B^*\bigr)
  \cdot\frac{1}{1+e^{\lambda_3(A^*-A_0)}}, \\[4pt]
\gamma_B B^* &= \bigl(g_B + g_{BA}A^*\bigr)
  \cdot\frac{1}{1+e^{\lambda_4(B^*-B_0)}}.
\end{align}
These coupled transcendental equations typically require numerical
solution.  Setting $\theta_B \approx g_{AB}B^*$ and
$\theta_A \approx g_{BA}A^*$ aligns the logistic inflection points
with the concentrations around which the system naturally operates,
which is particularly appropriate when transcription factors spend
most of their time near steady-state values.  Given steady-state-based
thresholds, the weighted-formulation matching yields
\begin{align}
\kappa_1 &= 2(g_A + \theta_B), &
\lambda_1 &= \frac{1}{\theta_B}
  \ln\!\left(1+\frac{2\theta_B}{g_A}\right), \\[4pt]
\kappa_2 &= 2(g_B + \theta_A), &
\lambda_2 &= \frac{1}{\theta_A}
  \ln\!\left(1+\frac{2\theta_A}{g_B}\right).
\end{align}
The choice between threshold selection strategies should be guided by
the availability of experimental data, the specific biological system
under study, and the intended modeling goals.

%----------------------------------------------------------------------
\subsection{Numerical Validation}
\label{subsubsec:numerical_val_params}
%----------------------------------------------------------------------

\begin{figure}[htbp]
\centering
\includegraphics[width=0.85\textwidth]{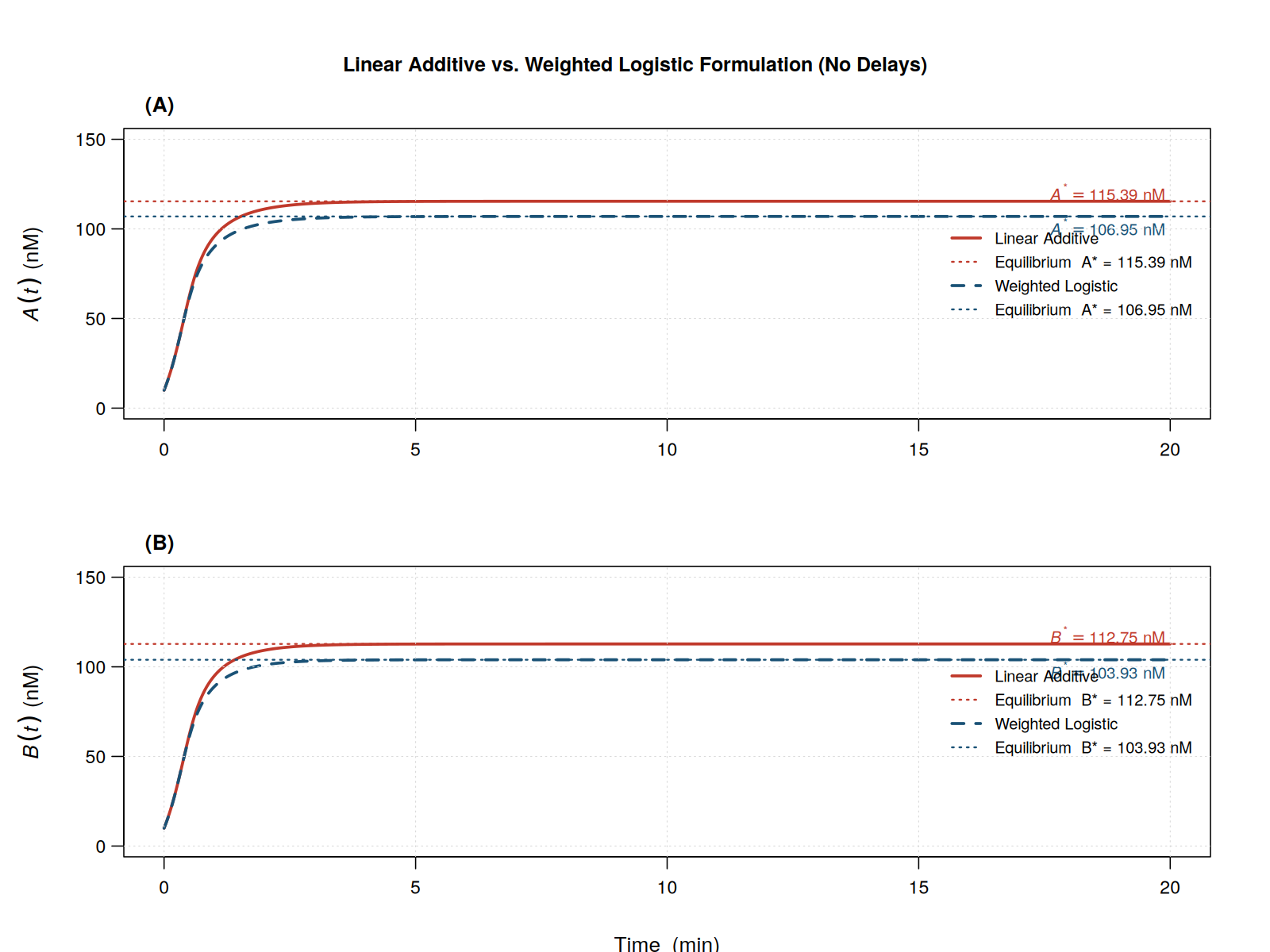}
\caption{Numerical comparison of the linear additive formulation
  (solid lines) and the weighted logistic reformulation (dashed lines)
  in the delay-free case, for protein~$A$ (top, red) and protein~$B$
  (bottom, blue). Both systems are integrated from the initial condition
  $(A(0), B(0)) = (10, 10)$~nM over 20~minutes.
  Parameters: $g_A = g_B = 50$~nM\,min$^{-1}$,
  $g_{AB} = g_{BA} = 2.5$~min$^{-1}$,
  $\gamma_A = 0.20$~min$^{-1}$, $\gamma_B = 0.24$~min$^{-1}$,
  repression thresholds $A_0 = B_0 = 70$~nM,
  logistic steepness $\lambda = n/A_0 = 4/70 \approx 0.057$~nM$^{-1}$;
  weighted logistic parameters $\kappa_1 = \kappa_2 = 200$~nM\,min$^{-1}$,
  $\lambda_1 = \lambda_2 = \ln 3/g_A \approx 0.022$~(nM\,min$^{-1}$)$^{-1}$,
  $\theta_A = \theta_B = g_A = 50$~nM\,min$^{-1}$. The two formulations
  converge to equilibria differing by approximately 7\%
  ($\approx 115$~nM vs.\ $\approx 107$~nM for protein~$A$), confirming
  that the $\approx 10\%$ slope overestimate identified in
  Remark~\ref{rem:approx_matching} has negligible impact on transient and
  steady-state dynamics. The weighted logistic trajectory converges
  slightly below the linear additive trajectory, consistent with
  saturation of the bounded activation term at high concentrations.
  Note: these illustrative parameters ($g_{AB} = g_{BA} = 2.5$~min$^{-1}$,
  $A_0 = B_0 = 70$~nM, $\lambda \approx 0.057$~nM$^{-1}$) differ from the
  main analysis parameters ($g_{AB} = g_{BA} = 3$~min$^{-1}$,
  $A_0 = B_0 = 100$~nM, $\lambda = 0.04$~nM$^{-1}$) used throughout
  Sections~\ref{sec:equilibrium}--\ref{sec:stability}; both sets lie within
  the experimentally relevant ranges of Table~\ref{tab:parameter_ranges}.}
\label{fig:logistic_compare_vinoth_2}
\end{figure}

Figure~\ref{fig:logistic_compare_vinoth_2} presents comparative
simulations of the linear additive formulation
(equations~\eqref{eq:A_explicit}--\eqref{eq:B_explicit} with logistic
self-repression) and the fully weighted logistic reformulation
(equations~\eqref{eq:A_explicit_basal}--\eqref{eq:B_explicit_basal}).
Using the parameter values
$g_A = g_B = 50~\mathrm{nM\,min^{-1}}$,
$g_{AB} = g_{BA} = 2.5~\mathrm{min^{-1}}$,
$\gamma_A = 0.20~\mathrm{min^{-1}}$,
$\gamma_B = 0.24~\mathrm{min^{-1}}$,
$A_0 = B_0 = 70~\mathrm{nM}$,
$\lambda_3 = \lambda_4 = n/A_0 \approx 0.057~\mathrm{nM^{-1}}$, and
initial conditions $A(0) = B(0) = 10~\mathrm{nM}$, both formulations
exhibit rapid, monotonic convergence to steady states.  These
parameters lie within the experimentally relevant ranges summarised
in Table~\ref{tab:parameter_ranges} and are consistent with values
reported for \textit{E.~coli} gene regulatory systems.  To facilitate
a clear comparison, attention is restricted here to the delay-free case
($\tau_1 = \tau_2 = \tau_{12} = \tau_{21} = 0$); the logistic DDE
extension is addressed in Section~\ref{sec:stability}.

The derived weighted logistic activation parameters
\begin{align*}
\kappa_1 &= \kappa_2 = 4g_A = 200~\mathrm{nM\,min^{-1}},
\quad
\lambda_1& = \lambda_2 = \frac{\ln 3}{g_A} \approx 0.022~\mathrm{(nM\,min^{-1})^{-1}}, \\
\theta_A &= \theta_B = g_A = 50~\mathrm{nM\,min^{-1}}
\end{align*}
ensure that the logistic approximation closely tracks the original
dynamics.  The agreement between the two trajectories confirms that
the parameter-matching strategy successfully preserves the system's
essential dynamical characteristics in moderate steepness regimes.

\begin{table}[htb!]
\small
\centering
\caption{Experimental ranges for parameters in the two-node gene
  regulatory network model with time delays.  Ranges are derived from
  synthetic GRNs in \textit{E.~coli}, reflecting variability in
  strains, inducers, and modifications
  \cite{vinoth2025extreme,madar2011negative,oehler1994quality}.}
\label{tab:parameter_ranges}
\begin{tabularx}{\textwidth}{l l X}
\toprule
\textbf{Parameter} & \textbf{Experimental Range} & \textbf{Description and Sources} \\
\midrule
$g_A$     & 10--100 nM\,min$^{-1}$  &
  Basal production rate for protein~$A$; from \textit{ara}/\textit{lac}
  basal expression \cite{madar2011negative,oehler1994quality}. \\
$g_B$     & 10--100 nM\,min$^{-1}$  &
  Symmetric to $g_A$; tuned for dynamics \cite{vinoth2025extreme}. \\
$g_{AB}$  & 1.0--5.0 min$^{-1}$     &
  Mutual activation strength; from AraC-mediated activation
  \cite{madar2011negative}. \\
$g_{BA}$  & 1.0--5.0 min$^{-1}$     &
  Symmetric activation; lac cooperative strengths
  \cite{oehler1994quality}. \\
$\gamma_A$& 0.05--0.30 min$^{-1}$   &
  Degradation rate; ssrA-tagged half-lives $\sim$2--14 min
  \cite{vinoth2025extreme}. \\
$\gamma_B$& 0.05--0.30 min$^{-1}$   &
  Symmetric degradation; consistent with synthetic oscillators. \\
$A_0$     & 50--200 nM              &
  Repression threshold; $K_d$ values from repressors
  \cite{madar2011negative}. \\
$B_0$     & 50--200 nM              &
  Symmetric threshold; lac $K_d$ variations \cite{oehler1994quality}. \\
$n$       & 2--5                    &
  Hill coefficient; effective cooperativity in \textit{ara}/\textit{lac}
  \cite{madar2011negative,oehler1994quality}. \\
$\lambda$ & 0.01--0.1 nM$^{-1}$    &
  Logistic steepness (fixed-weight formulation); derived from
  $n$ and thresholds via $\lambda = n/\theta$. \\
\bottomrule
\end{tabularx}
\end{table}

\begin{remark}[Delay-free restriction and the DDE extension]
\label{rem:delay_free}
Throughout this section we restrict attention to the delay-free case
($\tau_{12} = \tau_{21} = \tau_1 = \tau_2 = 0$) to facilitate a
clean comparison of equilibrium structure between the logistic
formulation and the original system.  The delay-dependent analysis
is taken up in Section~\ref{sec:stability}: local stability of both
equilibria is established analytically, and delay-induced Hopf
bifurcation is analysed for the symmetric configuration
$\tau_1 = \tau_2 = \tau_s$, yielding critical self-repression delays
$\tau_c \approx 1.19$~min (linear additive formulation) and
$\tau_c \approx 1.64$~min (weighted logistic formulation), together
with higher-order bifurcation sequences. 
\end{remark}

%%==========================================================================
\section{Equilibrium Analysis}
\label{sec:equilibrium}
%%==========================================================================

%----------------------------------------------------------------------
\subsection{Linear Additive Formulation}
\label{subsec:equil_linear}
%----------------------------------------------------------------------

The equilibrium points $(\bar{A}, \bar{B})$ of the linear additive
system~\eqref{eq:A_logistic}--\eqref{eq:B_logistic} satisfy
\begin{align}
0 &= \left( g_A + g_{AB} \bar{B} \right) f^-(\bar{A}, A_0, \lambda_A)
     - \gamma_A \bar{A}, \label{eq:eq_A}\\
0 &= \left( g_B + g_{BA} \bar{A} \right) f^-(\bar{B}, B_0, \lambda_B)
     - \gamma_B \bar{B}, \label{eq:eq_B}
\end{align}
which can be written in fixed-point form
\begin{align}
\bar{A} &= \frac{1}{\gamma_A} \left( g_A + g_{AB} \bar{B} \right)
            f^-(\bar{A}, A_0, \lambda_A), \label{eq:fixed_A}\\
\bar{B} &= \frac{1}{\gamma_B} \left( g_B + g_{BA} \bar{A} \right)
            f^-(\bar{B}, B_0, \lambda_B). \label{eq:fixed_B}
\end{align}
Since the right-hand side maps
$\Phi: \mathbb{R}^2_+ \to \mathbb{R}^2_+$ with bounded, globally smooth
logistic factors, the Brouwer fixed-point theorem guarantees at least one
equilibrium in any invariant rectangle. Uniqueness follows when the
operator norm $\|D\Phi\|_\infty < 1$.

Applying Newton--Raphson to~\eqref{eq:eq_A}--\eqref{eq:eq_B} with initial
guess $(A_0, B_0) = (100, 100)$~nM and convergence tolerance
$\varepsilon = 10^{-8}$ (implemented in R, package \texttt{nleqslv}), the
algorithm converged in 6 iterations to the unique biologically feasible
equilibrium
\begin{equation}
E^* = (\bar{A}, \bar{B}) = (167.96,\, 164.22)~\text{nM}.
\label{eq:equilibrium}
\end{equation}
Verification by direct substitution yields residuals below $10^{-10}$:
\begin{align}
|(g_A + g_{AB}\bar{B})f^-(\bar{A}, A_0, \lambda_A) - \gamma_A\bar{A}|
  &< 10^{-10}, \label{eq:residual1}\\
|(g_B + g_{BA}\bar{A})f^-(\bar{B}, B_0, \lambda_B) - \gamma_B\bar{B}|
  &< 10^{-10}. \label{eq:residual2}
\end{align}
The logistic function values at $E^*$ are:
\begin{align}
f_A^- &:= f^-(\bar{A}, A_0, \lambda_A)
       = \frac{1}{1+e^{0.04(167.96-100)}}
       = \frac{1}{1+e^{2.7184}}
       \approx 0.0619, \label{eq:f_minus_A}\\
f_B^- &:= f^-(\bar{B}, B_0, \lambda_B)
       = \frac{1}{1+e^{0.04(164.22-100)}}
       = \frac{1}{1+e^{2.5688}}
       \approx 0.0712. \label{eq:f_minus_B}
\end{align}
Both concentrations exceed the repression thresholds
($\bar{A} \approx 1.68\,A_0$, $\bar{B} \approx 1.64\,B_0$), so moderate
self-repression is active at equilibrium (6--7\% of maximal repression
function values).

%----------------------------------------------------------------------
\subsection{Weighted Logistic Formulation}
\label{subsec:equil_weighted}
%----------------------------------------------------------------------

Setting $\dot{A}=\dot{B}=0$ in~\eqref{eq:A_logistic_basal}--\eqref{eq:B_logistic_basal}
and suppressing all delays, the equilibrium
$E^*_{\mathrm{wl}}=(\bar{A}_{\mathrm{wl}},\bar{B}_{\mathrm{wl}})$ satisfies
\begin{align}
\bar{A} &= \frac{\kappa_1}{\gamma_A}
           f^+\!\bigl(g_{AB}\bar{B},\,\theta_B,\,\lambda_1\bigr)\,
           f^-\!\bigl(\bar{A},\,A_0,\,\lambda_3\bigr),
           \label{eq:wl_fixed_A}\\[4pt]
\bar{B} &= \frac{\kappa_2}{\gamma_B}
           f^+\!\bigl(g_{BA}\bar{A},\,\theta_A,\,\lambda_2\bigr)\,
           f^-\!\bigl(\bar{B},\,B_0,\,\lambda_4\bigr).
           \label{eq:wl_fixed_B}
\end{align}
Both right-hand sides are bounded (the product of two logistic functions
lies in $(0,1)$), so any equilibrium is confined to
$[0, \kappa_1/\gamma_A]\times[0, \kappa_2/\gamma_B] = [0,1000]\times[0,833]$~nM.
The Brouwer fixed-point theorem guarantees existence, and Newton--Raphson
with initial guess $(100,100)$~nM and tolerance $\varepsilon=10^{-10}$
yields the unique biologically feasible equilibrium
\begin{equation}
E^*_{\mathrm{wl}}
  = (\bar{A},\bar{B})
  = (144.46,\;139.99)~\text{nM}.
\label{eq:wl_equilibrium}
\end{equation}
The corresponding logistic function values are:
\begin{align}
f_1^+ &:= f^+\!\bigl(g_{AB}\bar{B},\theta_B,\lambda_1\bigr)
         = \frac{1}{1+e^{-0.02197(3\times139.99-50)}}
         = \frac{1}{1+e^{-8.129}}
         \approx 0.9997,
         \label{eq:f1plus_val}\\[4pt]
f_2^+ &:= f^+\!\bigl(g_{BA}\bar{A},\theta_A,\lambda_2\bigr)
         \approx 0.9998,
         \label{eq:f2plus_val}\\[4pt]
f_3^- &:= f^-\!\bigl(\bar{A},A_0,\lambda_3\bigr)
         = \frac{1}{1+e^{0.04(144.46-100)}}
         = \frac{1}{1+e^{1.7784}}
         \approx 0.1445,
         \label{eq:f3minus_val}\\[4pt]
f_4^- &:= f^-\!\bigl(\bar{B},B_0,\lambda_4\bigr)
         = \frac{1}{1+e^{0.04(139.99-100)}}
         = \frac{1}{1+e^{1.5996}}
         \approx 0.1680.
         \label{eq:f4minus_val}
\end{align}
The near-unity saturation of $f_1^+$ and $f_2^+$ (weighted signals
$g_{AB}\bar{B}\approx420$~nM\,min$^{-1}$ far exceed the activation threshold
$\theta_B=50$~nM\,min$^{-1}$) means that mutual activation is already in its
maximal regime at equilibrium.

\begin{remark}[Comparison of the two equilibria]
The linear additive formulation yields $E^*=(167.96,164.22)$~nM, roughly
16\% higher than $E^*_{\mathrm{wl}}=(144.46,139.99)$~nM. The shift arises
because the linear term $(g_A+g_{AB}B)$ grows without bound, whereas the
weighted logistic saturates at $\kappa_1=200$~nM\,min$^{-1}$; at
equilibrium the saturated logistic provides less activation, so
steady-state concentrations are lower. Both equilibria lie in the partial
self-repression regime, confirming that the qualitative regulatory
structure is preserved across formulations.
\end{remark}

%%==========================================================================
\section{Stability and Bifurcation Analysis}
\label{sec:stability}
%%==========================================================================

%----------------------------------------------------------------------
\subsection{Local Stability in the Delay-Free Case}
\label{subsec:stability_delay_free}
%----------------------------------------------------------------------

\noindent\textbf{Linear additive formulation.}\quad
Setting all delays to zero and linearizing~\eqref{eq:A_logistic}--\eqref{eq:B_logistic}
around $E^*$, the Jacobian is
\begin{equation}
J = \begin{pmatrix}
J_{11} & J_{12} \\
J_{21} & J_{22}
\end{pmatrix},
\label{eq:jacobian}
\end{equation}
with entries
\begin{align}
J_{11} &= (g_A + g_{AB} \bar{B})(-\lambda_A f_A^-(1-f_A^-)) - \gamma_A,
  \label{eq:J11}\\
J_{12} &= g_{AB} f_A^-, \label{eq:J12}\\
J_{21} &= g_{BA} f_B^-, \label{eq:J21}\\
J_{22} &= (g_B + g_{BA} \bar{A})(-\lambda_B f_B^-(1-f_B^-)) - \gamma_B,
  \label{eq:J22}
\end{align}
where $f_A^- := f^-(\bar{A},A_0,\lambda_A)$ and
$f_B^- := f^-(\bar{B},B_0,\lambda_B)$.
The characteristic polynomial is
\begin{equation}
\mu^2 - \mathrm{tr}(J)\mu + \det(J) = 0,
\label{eq:char_poly}
\end{equation}
with $\mathrm{tr}(J) = J_{11} + J_{22}$ and
$\det(J) = J_{11}J_{22} - J_{12}J_{21}$.
By the Routh--Hurwitz criterion, $E^*$ is locally asymptotically stable
if and only if $\mathrm{tr}(J) < 0$ and $\det(J) > 0$.

Expanding the trace:
\begin{align}
\mathrm{tr}(J) &= -\lambda_A(g_A+g_{AB}\bar{B})f_A^-(1-f_A^-)
                  -\lambda_B(g_B+g_{BA}\bar{A})f_B^-(1-f_B^-)
                  - (\gamma_A + \gamma_B). \label{eq:trace}
\end{align}
Since all parameters are strictly positive and $0 < f^\pm < 1$, every
term is negative, so $\mathrm{tr}(J) < 0$ holds \emph{unconditionally}
for any biologically relevant parameter choice.

The determinant is
\begin{align}
\det(J) &= \gamma_A\gamma_B
  + \lambda_B\gamma_A(g_B+g_{BA}\bar{A})f_B^-(1-f_B^-)
  + \lambda_A\gamma_B(g_A+g_{AB}\bar{B})f_A^-(1-f_A^-) \notag\\
&\quad
  + \lambda_A\lambda_B(g_A+g_{AB}\bar{B})(g_B+g_{BA}\bar{A})
    f_A^-(1-f_A^-)f_B^-(1-f_B^-)
  - g_{AB}g_{BA}f_A^-f_B^-. \label{eq:det_full}
\end{align}
Stability is maintained when mutual activation (the negative term) is
weak relative to self-repression and degradation.

Substituting the equilibrium values $E^* = (167.96, 164.22)$~nM and
$f_A^- \approx 0.0619$, $f_B^- \approx 0.0712$ at $E^*$:
\begin{align}
J_{11} &\approx -1.4605~\text{min}^{-1},   &  J_{12} &\approx 0.1857~\text{min}^{-1},\notag\\
J_{21} &\approx 0.2135~\text{min}^{-1},    &  J_{22} &\approx -1.7044~\text{min}^{-1}.
\label{eq:lin_J_numeric}
\end{align}
The trace, determinant, and discriminant are
\begin{align}
\mathrm{tr}(J) &\approx -3.1649~\text{min}^{-1} < 0, \notag\\
\det(J)        &\approx 2.4496~\text{min}^{-2} > 0, \notag\\
\Delta         &= \mathrm{tr}(J)^2 - 4\det(J) \approx 0.2181 > 0,
\label{eq:lin_traces}
\end{align}
yielding two real, negative eigenvalues
$\mu_1 \approx -1.349$~min$^{-1}$ and $\mu_2 \approx -1.816$~min$^{-1}$,
so $E^*$ is locally asymptotically stable in the delay-free case.

\begin{remark}[Absence of Hopf bifurcation in the delay-free case]
A Hopf bifurcation requires $\mathrm{tr}(J)=0$ with $\det(J)>0$.
Equation~\eqref{eq:trace} shows $\mathrm{tr}(J)<0$ universally, so
sustained oscillations require time delays, higher network dimensionality,
or additional nonlinear mechanisms.
\end{remark}

\noindent\textbf{Weighted logistic formulation.}\quad
Linearizing~\eqref{eq:A_logistic_basal}--\eqref{eq:B_logistic_basal}
around $E^*_{\mathrm{wl}}$ with all delays zero, the Jacobian is
$J^{\mathrm{wl}}$ with entries
\begin{align}
J_{11}^{\mathrm{wl}}
  &= -\kappa_1\lambda_3 f_1^+ f_3^-(1-f_3^-) - \gamma_A,
  \label{eq:wlJ11}\\[4pt]
J_{12}^{\mathrm{wl}}
  &= \kappa_1\lambda_1 g_{AB}\, f_1^+(1-f_1^+) f_3^-,
  \label{eq:wlJ12}\\[4pt]
J_{21}^{\mathrm{wl}}
  &= \kappa_2\lambda_2 g_{BA}\, f_2^+(1-f_2^+) f_4^-,
  \label{eq:wlJ21}\\[4pt]
J_{22}^{\mathrm{wl}}
  &= -\kappa_2\lambda_4 f_2^+ f_4^-(1-f_4^-) - \gamma_B.
  \label{eq:wlJ22}
\end{align}
Both diagonal entries are manifestly negative. The off-diagonal entries
are positive (mutual activation) but strongly attenuated at equilibrium
by the factor $f_i^+(1-f_i^+) \ll 1$ in the near-saturation regime.

The derivative products needed are:
\begin{align}
f_1^+(1-f_1^+) &\approx 2.95\times10^{-4},&
f_2^+(1-f_2^+) &\approx 2.20\times10^{-4},\notag\\
f_3^-(1-f_3^-) &\approx 0.1236,&
f_4^-(1-f_4^-) &\approx 0.1398.
\label{eq:deriv_prods}
\end{align}
Substituting numerically:
\begin{align}
J_{11}^{\mathrm{wl}}
  &= -200\times0.04\times0.9997\times0.1236 - 0.20
   \approx -1.1887~\text{min}^{-1},
   \label{eq:wlJ11_num}\\[4pt]
J_{12}^{\mathrm{wl}}
  &= 200\times0.02197\times3\times(2.95\times10^{-4})\times0.1445
   \approx 5.61\times10^{-4}~\text{min}^{-1},
   \label{eq:wlJ12_num}\\[4pt]
J_{21}^{\mathrm{wl}}
  &\approx 4.86\times10^{-4}~\text{min}^{-1},
   \label{eq:wlJ21_num}\\[4pt]
J_{22}^{\mathrm{wl}}
  &= -200\times0.04\times0.9998\times0.1398 - 0.24
   \approx -1.3581~\text{min}^{-1}.
   \label{eq:wlJ22_num}
\end{align}
The trace, determinant, and discriminant are:
\begin{align}
\mathrm{tr}(J^{\mathrm{wl}})
  &\approx -1.1887 + (-1.3581)
   = -2.5468~\text{min}^{-1} < 0,
   \label{eq:wl_trace}\\[4pt]
\det(J^{\mathrm{wl}})
  &\approx (-1.1887)(-1.3581) - (5.61\times10^{-4})(4.86\times10^{-4})
   \approx 1.6144~\text{min}^{-2} > 0,
   \label{eq:wl_det}\\[4pt]
\Delta &= \mathrm{tr}(J^{\mathrm{wl}})^2 - 4\det(J^{\mathrm{wl}})
       \approx 6.4861 - 6.4574 \approx 0.0287 > 0.
\label{eq:wl_disc}
\end{align}
Both eigenvalues are real and negative:
\begin{equation}
\mu_{1,2} = \frac{\mathrm{tr}(J^{\mathrm{wl}}) \pm \sqrt{\Delta}}{2}
  = \frac{-2.5468 \pm 0.1694}{2},
\label{eq:wl_eigenvalues}
\end{equation}
giving $\mu_1 \approx -1.1887$~min$^{-1}$ and $\mu_2 \approx -1.3581$~min$^{-1}$.
Because the off-diagonal coupling is $\mathcal{O}(10^{-4})$, the
eigenvalue shift relative to the diagonal entries is
$|J_{12}^{\mathrm{wl}}J_{21}^{\mathrm{wl}}|/|J_{11}^{\mathrm{wl}}-J_{22}^{\mathrm{wl}}|
\approx 1.6\times10^{-6}$~min$^{-1}$, so the eigenvalues coincide with the
diagonal entries $J_{11}^{\mathrm{wl}} \approx -1.1887$~min$^{-1}$ and
$J_{22}^{\mathrm{wl}} \approx -1.3581$~min$^{-1}$ to four-figure precision.
By the Routh--Hurwitz criterion, $E^*_{\mathrm{wl}}$ is \emph{locally
asymptotically stable} for all biologically relevant parameter values.

The near-zero off-diagonal coupling ($|J_{12}^{\mathrm{wl}}| \approx
5.6\times10^{-4}$~min$^{-1}$, two orders of magnitude below the diagonal
strengths $\approx 1.2$--$1.4$~min$^{-1}$) means that, near equilibrium,
the weighted logistic system behaves as two weakly coupled self-repressing
genes. Mutual activation becomes dynamically significant only during large
excursions from equilibrium --- precisely the regime responsible for the
extreme bursting events documented in~\cite{vinoth2025extreme}.

%----------------------------------------------------------------------
\subsection{Delay-Induced Hopf Bifurcation}
\label{subsec:hopf_analysis}
%----------------------------------------------------------------------

The introduction of nonzero time delays transforms the ODEs into DDEs with
transcendental characteristic equations, admitting Hopf bifurcations at
purely imaginary roots $\mu = i\omega$, $\omega > 0$.

\noindent\textbf{Linear additive formulation.}\quad
We focus on the symmetric delay configuration $\tau_1 = \tau_2 = \tau_s$,
$\tau_{12} = \tau_{21} = 0$.
The delay-dependent Jacobian entries at $E^*$ are
\begin{align}
J_{11}(i\omega) &= -\alpha_A e^{-i\omega\tau_s} - \gamma_A, &
J_{22}(i\omega) &= -\alpha_B e^{-i\omega\tau_s} - \gamma_B,
\end{align}
where the effective self-repression feedback strengths are
\begin{align}
\alpha_A &= \lambda_A(g_A + g_{AB}\bar{B})f_A^-(1-f_A^-)
         \approx 0.04\times542.66\times0.0581
         \approx 1.26~\text{min}^{-1}, \label{eq:alpha_A}\\
\alpha_B &= \lambda_B(g_B + g_{BA}\bar{A})f_B^-(1-f_B^-)
         \approx 0.04\times553.88\times0.0661
         \approx 1.46~\text{min}^{-1}, \label{eq:alpha_B}
\end{align}
and $J_{12}$, $J_{21}$ are real and delay-independent.
The mutual activation coupling is
\begin{equation}
\beta = g_{AB}g_{BA}f_A^-f_B^-
      = 9\times0.0619\times0.0712
      \approx 0.0397. \label{eq:beta}
\end{equation}

Setting $\mu = i\omega$ in the characteristic equation
$\det(\mu I - J(\mu))=0$ gives
\begin{equation}
\bigl(i\omega + \gamma_A + \alpha_A e^{-i\omega\tau_s}\bigr)
\bigl(i\omega + \gamma_B + \alpha_B e^{-i\omega\tau_s}\bigr)
- \beta = 0.
\label{eq:char_eq_iomega}
\end{equation}
Writing $c = \cos(\omega\tau_s)$, $s = \sin(\omega\tau_s)$ and
$P_k = (\gamma_k + \alpha_k c) + i(\omega - \alpha_k s)$ for $k\in\{A,B\}$,
equation~\eqref{eq:char_eq_iomega} becomes $P_A P_B = \beta$.
Separating real and imaginary parts yields the \emph{two-equation
transcendental system} that must be satisfied simultaneously:
\begin{subequations}\label{eq:hopf_system}
\begin{align}
(\gamma_A + \alpha_A c)(\gamma_B + \alpha_B c)
  - (\omega - \alpha_A s)(\omega - \alpha_B s)
  &= \beta,
\label{eq:hopf_re}\\
\omega\bigl[(\gamma_A+\gamma_B)+(\alpha_A+\alpha_B)c\bigr]
  - s\bigl[\alpha_A\gamma_B + \alpha_B\gamma_A + 2\alpha_A\alpha_B c\bigr]
  &= 0.
\label{eq:hopf_im}
\end{align}
\end{subequations}

\begin{remark}[The imaginary-trace shortcut fails for asymmetric parameters]
\label{rem:trace_failure}
One might attempt to simplify~\eqref{eq:hopf_im} by setting $s = 0$
(i.e.\ $\omega\tau_s = n\pi$). With $s=0$ and $c=\pm 1$,
equation~\eqref{eq:hopf_im} reduces to
$\omega[(\gamma_A+\gamma_B)\pm(\alpha_A+\alpha_B)]=0$, which requires
$(\gamma_A+\gamma_B) = \mp(\alpha_A+\alpha_B)$, i.e.\ $0.44 = \mp 2.72$.
This is violated by more than a factor of six for the present parameters,
so $\omega\tau_s = n\pi$ is \emph{not} a solution of the imaginary-part
equation~\eqref{eq:hopf_im}, and the two equations must be solved jointly.
\end{remark}

System~\eqref{eq:hopf_system} is solved numerically in R using the
\texttt{nleqslv} package, initialised on a fine grid of
$(\omega,\,\theta)$ with $\theta=\omega\tau_s\in(0,\pi)$,
and residuals verified against the original complex
equation~\eqref{eq:char_eq_iomega} to machine precision.
The smallest positive $\tau_s$ satisfying~\eqref{eq:hopf_system}
is the first Hopf bifurcation point:
\begin{equation}
\omega_c \approx 1.3519~\text{min}^{-1}, \qquad
\tau_c   \approx 1.1879~\text{min} \approx 71.3~\text{s}, \qquad
\theta_c = \omega_c\tau_c \approx 1.6059~\text{rad} \approx 0.5112\pi.
\label{eq:critical_delay}
\end{equation}
A direct substitution into~\eqref{eq:char_eq_iomega}
yields $|\mathrm{Re}|<10^{-15}$, $|\mathrm{Im}|<10^{-15}$.

The transversality condition,
evaluated via the implicit function theorem
($d\mu/d\tau = -(\partial F/\partial\tau)/(\partial F/\partial\mu)$
with $F=P_AP_B-\beta$),
\begin{equation}
\left.\frac{d\operatorname{Re}(\mu)}{d\tau}\right|_{\tau=\tau_c}
\approx 0.5150~\text{min}^{-2} \neq 0,
\label{eq:transversality_lin}
\end{equation}
confirms a non-degenerate Hopf bifurcation.
The oscillation period at onset is
\begin{equation}
T = \frac{2\pi}{\omega_c}
  \approx \frac{6.2832}{1.3519}
  \approx 4.648~\text{min}.
\label{eq:period_lin}
\end{equation}

\begin{remark}[Biological context of $\tau_c$]
The critical delay $\tau_c \approx 1.19$~min is consistent with the
faster end of transcription-translation lags in bacterial and fast-cycling
eukaryotic networks: mRNA synthesis of short genes (0.3--0.8~min),
translation of small proteins (0.3--0.8~min), and nuclear export or
folding ($\leq 0.5$~min) can collectively produce delays of this
magnitude~\cite{madar2011negative,oehler1994quality}.
\end{remark}

\noindent\textbf{Weighted logistic formulation.}\quad
For the weighted logistic system with $\tau_1=\tau_2=\tau_s>0$ and
$\tau_{12}=\tau_{21}=0$, the effective self-repression feedback strengths
at $E^*_{\mathrm{wl}}$ are
\begin{align}
\alpha_A^{\mathrm{wl}}
  &:= \kappa_1\lambda_3 f_1^+ f_3^-(1-f_3^-)
   = 200\times0.04\times0.9997\times0.1236
   \approx 0.9887~\text{min}^{-1},
   \label{eq:wl_alphaA}\\[4pt]
\alpha_B^{\mathrm{wl}}
  &:= \kappa_2\lambda_4 f_2^+ f_4^-(1-f_4^-)
   = 200\times0.04\times0.9998\times0.1398
   \approx 1.1181~\text{min}^{-1}.
   \label{eq:wl_alphaB}
\end{align}
The off-diagonal coupling term is
$J_{12}^{\mathrm{wl}}J_{21}^{\mathrm{wl}} \approx 2.73\times10^{-7}$,
negligible in all computations below.

The characteristic equation at $\mu=i\omega$ takes the same form
as~\eqref{eq:char_eq_iomega} with $\alpha_k$ replaced by
$\alpha_k^{\mathrm{wl}}$ and $\beta$ replaced by
$J_{12}^{\mathrm{wl}}J_{21}^{\mathrm{wl}}\approx 2.73\times10^{-7}$
(hereafter set to zero without loss of precision).
Writing $c=\cos(\omega\tau_s)$, $s=\sin(\omega\tau_s)$, the
two-equation transcendental system analogous
to~\eqref{eq:hopf_system} is
\begin{subequations}\label{eq:hopf_system_wl}
\begin{align}
\bigl(\gamma_A + \alpha_A^{\mathrm{wl}}c\bigr)
\bigl(\gamma_B + \alpha_B^{\mathrm{wl}}c\bigr)
- \bigl(\omega - \alpha_A^{\mathrm{wl}}s\bigr)
  \bigl(\omega - \alpha_B^{\mathrm{wl}}s\bigr)
&= 0,
\label{eq:hopf_wl_re}\\
\omega\bigl[(\gamma_A+\gamma_B)
  +(\alpha_A^{\mathrm{wl}}+\alpha_B^{\mathrm{wl}})c\bigr]
- s\bigl[\alpha_A^{\mathrm{wl}}\gamma_B
       + \alpha_B^{\mathrm{wl}}\gamma_A
       + 2\alpha_A^{\mathrm{wl}}\alpha_B^{\mathrm{wl}}c\bigr]
&= 0.
\label{eq:hopf_wl_im}
\end{align}
\end{subequations}

\begin{remark}[Imaginary-trace shortcut fails here too]
\label{rem:trace_failure_wl}
Setting $s=0$ in~\eqref{eq:hopf_wl_im} requires
$(\gamma_A+\gamma_B)\pm(\alpha_A^{\mathrm{wl}}+\alpha_B^{\mathrm{wl}})=0$,
i.e.\ $0.44 = \mp 2.107$, violated by a factor of five.
The two equations must be solved jointly.
\end{remark}

System~\eqref{eq:hopf_system_wl} is solved numerically with the same
grid-initialised Newton procedure used for the linear additive case,
and residuals are verified against the complex characteristic equation to
machine precision.
The smallest positive $\tau_s$ satisfying~\eqref{eq:hopf_system_wl}
is the first Hopf bifurcation point:
\begin{equation}
\omega_c^{\mathrm{wl}} \approx 1.0920~\text{min}^{-1}, \quad
\tau_c^{\mathrm{wl}}   \approx 1.637~\text{min} \approx 98.2~\text{s}, \quad
\theta_c^{\mathrm{wl}} = \omega_c^{\mathrm{wl}}\tau_c^{\mathrm{wl}}
  \approx 1.787~\text{rad} \approx 0.569\pi.
\label{eq:wl_critical_delay}
\end{equation}
A direct substitution into the characteristic equation
yields $|\mathrm{Re}|<10^{-11}$, $|\mathrm{Im}|<10^{-11}$.

The transversality condition, evaluated via the implicit function
theorem ($d\mu/d\tau=-(\partial F/\partial\tau)/(\partial F/\partial\mu)$),
\begin{equation}
\left.\frac{d\operatorname{Re}(\mu)}{d\tau}\right|_{\tau=\tau_c^{\mathrm{wl}}}
\approx 0.232~\text{min}^{-2} \neq 0,
\label{eq:transversality_wl}
\end{equation}
confirms a non-degenerate Hopf bifurcation.
The oscillation period at onset is
\begin{equation}
T^{\mathrm{wl}} = \frac{2\pi}{\omega_c^{\mathrm{wl}}}
  \approx \frac{6.2832}{1.0920}
  \approx 5.754~\text{min}.
\label{eq:wl_period}
\end{equation}

\begin{remark}[Critical delay comparison across formulations]
\label{rem:tau_c_comparison}
The weighted logistic formulation yields
$\tau_c^{\mathrm{wl}}\approx1.64$~min, a $\approx37\%$ increase over
$\tau_c\approx1.19$~min for the linear additive formulation. The increase
arises from the smaller effective self-repression feedback strengths
($\alpha_A^{\mathrm{wl}}\approx0.99$ vs.\
$\alpha_A^{\mathrm{lin}}\approx1.26$~min$^{-1}$), which in turn reflect
the lower equilibrium concentrations: at $\bar{A}=144.46$~nM, the
self-repression logistic is less deeply repressed
($f_3^-\approx 0.1445$) than at $\bar{A}=167.96$~nM
($f_A^-\approx 0.0619$) in the linear additive case, but the bounded
weighted-logistic activation produces a much smaller pre-factor
$\kappa_1 f_1^+\approx 200$~nM\,min$^{-1}$ than the unbounded linear
combination $g_A+g_{AB}\bar{B}\approx 542.66$~nM\,min$^{-1}$, so the
net feedback strength $\alpha=\lambda(\text{pre-factor})f^-(1-f^-)$ is
smaller. Both critical delays fall within the biologically plausible
range of transcription-translation lags.
\end{remark}

%----------------------------------------------------------------------
\subsection{Higher-Order Bifurcations}
\label{subsec:higher_order}
%----------------------------------------------------------------------

Once the primary critical pair $(\omega_c,\tau_c)$ has been identified,
the structure of the transcendental characteristic
equation~\eqref{eq:char_eq_iomega} immediately yields a countable family
of additional Hopf points. If $\mu = i\omega_c$ is a root for delay
$\tau = \tau_c$, then it remains a root for every delay
\begin{equation}
\tau_c^{(k)} = \tau_c + \frac{2k\pi}{\omega_c},
\qquad k = 0, 1, 2, \dots,
\label{eq:tau_c_higher}
\end{equation}
because $e^{-i\omega_c \tau_c^{(k)}} = e^{-i\omega_c\tau_c}e^{-2k\pi i}
= e^{-i\omega_c\tau_c}$.  Each successive value $\tau_c^{(k)}$ corresponds
to the same imaginary root re-entering the characteristic equation as the
delay is increased, and therefore lies on the branch
$\theta = \omega_c\tau_c + 2k\pi$ in the
$(\omega,\theta)$ formulation~\eqref{eq:hopf_system}.

\noindent\textbf{Linear additive formulation.}\quad
With $(\omega_c,\tau_c) \approx (1.3519~\text{min}^{-1},\,1.188~\text{min})$,
formula~\eqref{eq:tau_c_higher} gives
\begin{equation}
\tau_c^{(1)} \approx 1.188 + \frac{2\pi}{1.3519}
  \approx 5.836~\text{min}, \qquad
\tau_c^{(2)} \approx 10.483~\text{min}, \qquad
\tau_c^{(3)} \approx 15.131~\text{min}.
\label{eq:lin_higher_taus}
\end{equation}

A second sequence of Hopf points arises from a
\emph{distinct} imaginary root of~\eqref{eq:hopf_system}.  Numerical
exploration of the same branch reveals a secondary solution at
\begin{equation}
\omega_c' \approx 1.2989~\text{min}^{-1}, \qquad
\tau_c' \approx 1.4351~\text{min},
\label{eq:lin_secondary}
\end{equation}
whose own family of replicas
$\tau_c'{}^{(k)} = \tau_c' + 2k\pi/\omega_c'$ interleaves with the primary
sequence.

\noindent\textbf{Weighted logistic formulation.}\quad
With $(\omega_c^{\mathrm{wl}},\tau_c^{\mathrm{wl}})
\approx (1.0920,\,1.637)$~min, the analogous higher Hopf points are
\begin{equation}
\tau_c^{\mathrm{wl},(1)} \approx 1.637 + \frac{2\pi}{1.0920}
  \approx 7.39~\text{min}, \qquad
\tau_c^{\mathrm{wl},(2)} \approx 13.14~\text{min}, \qquad
\tau_c^{\mathrm{wl},(3)} \approx 18.90~\text{min},
\label{eq:wl_higher_taus}
\end{equation}
with a secondary sequence seeded at
$(\omega_c'{}^{\mathrm{wl}},\,\tau_c'{}^{\mathrm{wl}})
\approx (0.9682,\,1.833)$~min.
Both higher-order sequences~\eqref{eq:lin_higher_taus}
and~\eqref{eq:wl_higher_taus} have been verified numerically by direct
substitution into the characteristic equation~\eqref{eq:char_eq_iomega}
to machine precision.

\begin{remark}[Stability between Hopf points]
Between consecutive Hopf points each crossing of the imaginary axis
adds a pair of unstable eigenvalues (a standard consequence of the
positive transversality condition).  Stability is therefore not
recovered for $\tau > \tau_c$: the equilibrium remains unstable for all
$\tau \geq \tau_c$, and the higher-order
points~\eqref{eq:lin_higher_taus}--\eqref{eq:wl_higher_taus} mark the
appearance of additional unstable directions rather than a return to
stability.  This monotone destabilisation is characteristic of negative
delayed feedback and underlies the period-doubling and chaotic regimes
documented at large delays in~\cite{vinoth2025extreme}.
\end{remark}

%----------------------------------------------------------------------
\subsection{Biological Implications of the Delay-Induced Hopf Bifurcation}
\label{subsec:bio_implications}
%----------------------------------------------------------------------

The critical delays $\tau_c \approx 1.19$~min (linear additive) and
$\tau_c^{\mathrm{wl}} \approx 1.64$~min (weighted logistic) carry several
important consequences. First, the combination of strong mutual activation
with self-repression creates an intrinsic oscillatory tendency requiring
delays of only one to two minutes to destabilize the equilibrium, explaining
why oscillatory behavior is ubiquitous in gene networks where regulatory
steps accumulate temporal lags. Second, for synthetic biology applications,
maintaining stability in engineered circuits requires either minimizing
transcription-translation delays through fast-degrading mRNA species,
reducing mutual activation strength, or increasing self-repression steepness
through larger $\lambda$. Third, delays exceeding the critical value produce
sustained oscillations with period $T = 2\pi/\omega_c \approx 4.65$~min for
the linear additive formulation and
$T^{\mathrm{wl}} = 2\pi/\omega_c^{\mathrm{wl}} \approx 5.75$~min for the
weighted logistic formulation; further increases are conjectured to drive
period-doubling cascades toward chaos and extreme bursting events analogous
to those documented for the original Hill-based formulation
in~\cite{vinoth2025extreme}, although direct numerical confirmation in the
logistic reformulations constitutes important future work.

A critical observation undergirds the entire analysis: all Hopf bifurcation
threshold values are determined by the local Jacobian structure at the
equilibrium. Because the parameter matching
$\lambda_A = \lambda_B = n/\theta_0 = 0.04$~nM$^{-1}$ ensures equivalent
Jacobian eigenvalues at the equilibrium, the \emph{onset} of oscillatory
instability is preserved under the Hill-to-logistic transformation. The rich
complex dynamics documented in~\cite{vinoth2025extreme} at larger delay
values --- period-doubling cascades, intermittency, and crisis-induced
extreme events --- arise from global nonlinear dynamics far from equilibrium,
and their preservation under the logistic reformulation requires direct
numerical validation, which constitutes important future work.

%%==========================================================================
\section{Lipschitz Constants and Numerical Efficiency}
\label{sec:lipschitz}
%%==========================================================================

Global Lipschitz constants bound trajectory divergence and determine
integration step-size requirements for numerical solvers. We compute them
explicitly for both formulations, beginning with concentration bounds for the
biologically relevant domain. At the linear additive equilibrium
$E^* = (167.96, 164.22)$~nM, both proteins are approximately 1.7 times
their repression thresholds and the logistic repression functions have
decreased to 6--7\% of their maximum values.

Allowing for a 150\% overshoot from equilibrium and noting that at
concentrations significantly above $1.7\theta$ the repression term $f^-$
provides exponentially strong suppression (e.g.\
$f^-(4.2\,\theta) \approx 2.8\times10^{-6}$), we adopt the conservative
estimate
\begin{equation}
A_{\max} = B_{\max} = 500~\text{nM}.
\label{eq:max_concentrations}
\end{equation}
This bound is biologically justified and valid across a wide range of
parameter perturbations.

%----------------------------------------------------------------------
\subsection{Linear Additive Formulation}
\label{subsec:lipschitz_linear}
%----------------------------------------------------------------------

The Jacobian of system~\eqref{eq:A_logistic}--\eqref{eq:B_logistic}, used to bound the vector-field Lipschitz constant $L_F$, has entries
\begin{align}
J_{11} &= (g_A + g_{AB}B)(-\lambda_A f_A^-(1-f_A^-)) - \gamma_A, &
J_{12} &= g_{AB} f_A^-, \notag\\
J_{21} &= g_{BA} f_B^-, &
J_{22} &= (g_B + g_{BA}A)(-\lambda_B f_B^-(1-f_B^-)) - \gamma_B.
\end{align}
Using $f^-(1-f^-)\leq 1/4$ and $0<f^-<1$:
\begin{align}
|J_{11}| &\leq (g_A + g_{AB}B_{\max})\cdot\frac{\lambda_A}{4} + \gamma_A
           = (50 + 3\times500)\times0.01 + 0.20
           = 15.70~\text{min}^{-1}, \notag\\
|J_{12}| &\leq g_{AB} = 3.00~\text{min}^{-1}, \notag\\
|J_{21}| &\leq g_{BA} = 3.00~\text{min}^{-1}, \notag\\
|J_{22}| &\leq (g_B + g_{BA}A_{\max})\cdot\frac{\lambda_B}{4} + \gamma_B
           = 1550\times0.01 + 0.24
           = 15.74~\text{min}^{-1}.
\end{align}
The global Lipschitz constant (infinity norm):
\begin{align}
L_F^{\mathrm{lin}}
  &\leq \max\!\bigl(|J_{11}|+|J_{12}|,\;|J_{21}|+|J_{22}|\bigr)
   = \max(15.70+3.00,\;3.00+15.74)
   = \boxed{18.74~\text{min}^{-1}}.
\label{eq:lipschitz_F}
\end{align}

Turning to the Jacobian Lipschitz constant $L_{DF}$, the second derivative of the decreasing logistic satisfies
\begin{equation}
\left|\frac{\partial^2 f^-}{\partial x^2}\right|
  \leq \lambda^2\rho,
  \quad
  \rho = \max_{f\in(0,1)}|f(1-f)(1-2f)|
       = \frac{\sqrt{3}}{18} \approx 0.09623.
\end{equation}
The dominant second-order partial is
\begin{align}
\left|\frac{\partial^2 \dot{A}}{\partial A^2}\right|
  &\leq (g_A + g_{AB}B_{\max})\cdot\lambda_A^2\cdot\rho
   = 1550\times(0.04)^2\times0.09623
   \approx 0.238~\text{nM}^{-1}\text{min}^{-1},
\end{align}
and the cross-partial satisfies
\begin{align}
\left|\frac{\partial^2 \dot{A}}{\partial A\,\partial B}\right|
  &\leq g_{AB}\cdot\frac{\lambda_A}{4}
   = 3\times0.01
   = 0.030~\text{nM}^{-1}\text{min}^{-1}.
\end{align}
The Jacobian Lipschitz constant:
\begin{equation}
L_{DF}^{\mathrm{lin}} = \max\{0.238, 0.030\}
  = \boxed{0.24~\text{nM}^{-1}\text{min}^{-1}}.
\label{eq:lipschitz_DF}
\end{equation}

\begin{remark}[Comparison with Hill formulations]
For Hill functions with non-integer exponent $n\in(1,2)$, the second
derivative contains terms proportional to $x^{n-2}$, which diverge as
$x\to0^+$:
$\partial^2 h^+/\partial x^2 \sim (n(n-1)/\theta^n)\,x^{n-2} \to \infty$.
This singularity forces $L_{DF}^{\text{Hill}} = \infty$ on any domain
containing the origin, requiring special adaptive integration strategies
in low-expression regimes. The logistic formulation's finite
$L_{DF}^{\mathrm{lin}} = 0.24$ ensures uniform well-conditioning
throughout the entire biologically relevant domain.
\end{remark}

%----------------------------------------------------------------------
\subsection{Weighted Logistic Formulation}
\label{subsec:lipschitz_weighted}
%----------------------------------------------------------------------

The Jacobian of system~\eqref{eq:A_logistic_basal}--\eqref{eq:B_logistic_basal}, used to bound $L_F$, has entries
\begin{align}
J_{11}^{\mathrm{wt}} &= -\kappa_1\lambda_3 f_1^+ f_3^-(1-f_3^-)-\gamma_A, &
J_{12}^{\mathrm{wt}} &= \kappa_1\lambda_1 g_{AB} f_1^+(1-f_1^+)f_3^-, \notag\\
J_{21}^{\mathrm{wt}} &= \kappa_2\lambda_2 g_{BA} f_2^+(1-f_2^+)f_4^-, &
J_{22}^{\mathrm{wt}} &= -\kappa_2\lambda_4 f_2^+ f_4^-(1-f_4^-)-\gamma_B.
\end{align}
Bounds:
\begin{align}
|J_{11}^{\mathrm{wt}}| &\leq \kappa_1\lambda_3\cdot\tfrac{1}{4} + \gamma_A
  = 200\times0.04\times0.25 + 0.20
  = 2.20~\text{min}^{-1}, \notag\\
|J_{12}^{\mathrm{wt}}| &\leq \kappa_1\lambda_1 g_{AB}\cdot\tfrac{1}{4}
  = 200\times0.02197\times3\times0.25
  \approx 3.30~\text{min}^{-1}, \notag\\
|J_{21}^{\mathrm{wt}}| &\approx 3.30~\text{min}^{-1}, \notag\\
|J_{22}^{\mathrm{wt}}| &\leq \kappa_2\lambda_4\cdot\tfrac{1}{4} + \gamma_B
  = 200\times0.04\times0.25 + 0.24
  = 2.24~\text{min}^{-1}.
\end{align}
\begin{equation}
L_F^{\mathrm{wt}}
  \leq \max(2.20+3.30,\;3.30+2.24)
  = \boxed{5.54~\text{min}^{-1}},
\label{eq:lipschitz_F_wt}
\end{equation}
an approximately 70\% reduction relative to $L_F^{\mathrm{lin}}=18.74$.

For the Jacobian Lipschitz constant $L_{DF}$, the dominant second-order partial for the weighted logistic is:
\begin{align}
\left|\frac{\partial^2 \dot{A}}{\partial B^2}\right|
  &\leq \kappa_1(\lambda_1 g_{AB})^2\rho
   = 200\times(0.02197\times3)^2\times0.09623
   \approx 0.0836~\text{nM}^{-1}\text{min}^{-1}, \notag\\
\left|\frac{\partial^2 \dot{A}}{\partial A^2}\right|
  &\leq \kappa_1\lambda_3^2\rho
   = 200\times(0.04)^2\times0.09623
   \approx 0.0308~\text{nM}^{-1}\text{min}^{-1}, \notag\\
\left|\frac{\partial^2 \dot{A}}{\partial A\,\partial B}\right|
  &\leq \kappa_1\lambda_1 g_{AB}\cdot\tfrac{1}{4}\cdot\lambda_3\cdot\tfrac{1}{4}
   \approx 0.033~\text{nM}^{-1}\text{min}^{-1}.
\end{align}
\begin{equation}
L_{DF}^{\mathrm{wt}}
  = \max\{0.0836, 0.0308, 0.033\}
  = \boxed{0.084~\text{nM}^{-1}\text{min}^{-1}},
\label{eq:lipschitz_DF_wt}
\end{equation}
a 65\% reduction relative to $L_{DF}^{\mathrm{lin}}=0.24$.

%----------------------------------------------------------------------
\subsection{Comparative Summary}
\label{subsec:lipschitz_comparison}
%----------------------------------------------------------------------

Table~\ref{tab:lipschitz_comparison} summarises the quantitative comparison.

\begin{table}[h]
\centering
\caption{Lipschitz constants: linear additive vs.\ weighted logistic
  formulations. $L_{DF}$ values are rounded to 2 significant figures;
  exact computed values are $0.238$ and $0.0836$~nM$^{-1}$\,min$^{-1}$
  respectively.}
\label{tab:lipschitz_comparison}
\begin{tabular}{lccc}
\toprule
\textbf{Constant}
  & \textbf{Linear Additive}
  & \textbf{Weighted Logistic}
  & \textbf{Ratio (Wt/Lin)} \\
\midrule
$L_F$ (min$^{-1}$)              & 18.74 & 5.54  & 0.30 \\
$L_{DF}$ (nM$^{-1}$ min$^{-1}$) & 0.24  & 0.084 & 0.35 \\
$L_{DF}/L_F$ (nM$^{-1}$)        & 0.013 & 0.015 & 1.18 \\
\bottomrule
\end{tabular}
\end{table}

The approximately threefold reduction in both Lipschitz constants stems
from replacing unbounded linear activation with bounded sigmoidal activation:
as $g_{AB}B\to\infty$, $f_1^+\to1$, preventing runaway growth and bounding
all derivatives. The saturation reflects real cellular constraints ---
finite RNA polymerase availability, ribosome competition, and metabolic
burden --- that the linear additive model ignores.

From a computational perspective, Runge--Kutta methods require step sizes
inversely proportional to the Lipschitz constant for stability; the weighted
logistic formulation permits steps approximately three times larger for
comparable accuracy. The chain of bounds
\begin{equation}
\|\mathbf{F}\|_\infty \leq M, \quad
\|\nabla\mathbf{F}\|_\infty \leq L_F, \quad
\|\nabla^2\mathbf{F}\|_\infty \leq L_{DF}
\end{equation}
establishes that both vector fields are globally $C^2$ with uniformly bounded
derivatives, ensuring that standard existence/uniqueness theorems apply and
that numerical methods converge at their theoretical rates throughout the
biologically relevant domain.

The ratio $L_{DF}/L_F$ is nearly identical between formulations
($0.013$ vs.\ $0.015$~nM$^{-1}$), indicating comparable relative
sensitivity structure. What differs is the absolute magnitude: the weighted
logistic system produces proportionally smaller responses to state
perturbations across the entire concentration domain.

The nonzero basal expression in the activator-free regime,
\begin{equation}
f^+(0,\theta,\lambda) = \frac{1}{1 + e^{\lambda\theta}} > 0,
\label{eq:basal_expression}
\end{equation}
guarantees continuous control authority even during pre-burst low-expression
phases, enabling predictive intervention strategies that can anticipate and
prevent extreme events before they fully manifest. This eliminates the
controllability gaps that Hill functions create at zero concentration, where
zero production precludes any recovery through intrinsic dynamics alone.

%%==========================================================================
\section{Conclusion}
%%==========================================================================

This paper has argued, through systematic theoretical analysis of a
delay-coupled two-gene regulatory network studied under two distinct logistic
reformulations, that logistic functions constitute a mathematically principled,
computationally efficient, and biologically realistic alternative to Hill
functions for modeling gene regulatory networks. The hybrid regulatory
architecture of Vinoth et al.~\cite{vinoth2025extreme}---a linear activation
$(g_A + g_{AB}B)$ tracing its conceptual origins to Jacob and Monod's seminal
operon model~\cite{jacob1961genetic}, combined with Hill self-repression---admits
a rigorous transformation into a purely sigmoidal framework by supplanting the
unbounded linear term with a bounded logistic function engineered to reproduce
both the desired basal rate and the initial activation slope via Taylor expansion
matching. The advantages of this substitution are not incidental but arise from
deep structural properties of the logistic form, manifesting simultaneously
across the mathematical, numerical, and biological dimensions that matter most
in practice.

From a mathematical standpoint, the global infinite differentiability of
logistic functions removes the singularities that afflict Hill functions
whenever the cooperativity exponent takes a non-integer value, a ubiquitous
occurrence when fitting experimental dose-response data. The self-referential
derivative structure $f' = \lambda f(1-f)$ significantly simplifies Jacobian
computations for stability analysis, reducing Jacobian entries to products of
already-evaluated function values and eliminating fractional exponents
entirely. The full decoupling of threshold $\theta$ and steepness
$\lambda$ allows each to be tuned independently, a practical advantage that
Hill functions, where both properties are entangled through the coefficient
$n/(4\theta)$, cannot provide. These structural properties translate directly
into quantitative numerical gains: replacing linear additive activation by
weighted logistic activation in the Vinoth et al.\ two-gene delay network
reduces the global Lipschitz constant from $L_F = 18.74$~min$^{-1}$ to
$5.54$~min$^{-1}$, an approximately 70\% reduction, and the Jacobian Lipschitz
constant from $L_{DF} = 0.24$ to $0.084$~nM$^{-1}$\,min$^{-1}$, a 65\%
reduction. Both bounds are finite and uniform across the entire biologically
relevant state space, in contrast to Hill formulations with non-integer
exponents where $L_{DF}$ diverges as trajectories approach zero.

The core case study analyzed here, the mutual-activation and self-repression
delay-coupled two-gene network, illustrates the framework's analytical power
through two parallel logistic reformulations. The parameter matching
$\lambda = n/\theta$, derived by equating logistic and Hill slopes at their
half-maximal points, ensures local equivalence to the original Hill-based model
at each respective equilibrium. For the linear additive logistic formulation,
the unique biologically feasible equilibrium $(167.96,\,164.22)$~nM was computed
numerically and verified analytically; the weighted logistic formulation yields
the distinct equilibrium $(144.46,\,139.99)$~nM, which is lower because the
bounded activation saturates at high concentrations rather than growing without
bound. Local stability analysis established that the Jacobian trace is strictly
negative for all positive parameters in both formulations in the delay-free
case ($\tau=0$); with delays, both formulations remain stable for
$\tau\in(0,\tau_c)$ and lose stability via Hopf bifurcation at $\tau_c$.
Hopf bifurcation analysis identified critical self-repression delays of
$\tau_c \approx 1.19$~min for the linear additive formulation and
$\tau_c \approx 1.64$~min for the weighted logistic formulation, with
higher-order bifurcation sequences characterized analytically in each case.
Both critical delays fall within the biologically plausible range of
transcription-translation lags. Direct numerical simulation of complex dynamics
and extreme events in the logistic reformulations remains an important direction
for future work.

Basal expression, structurally absent from Hill functions and typically patched
through ad~hoc additive offsets that lack biological justification and distort
the normalized response range, emerges as an intrinsic property of the logistic
form. Its magnitude is governed by existing parameters through the strictly
positive value $f^+(0,\theta,\lambda) = 1/(1+e^{\lambda\theta})$, requiring no
additional parameters. The decreasing logistic function for repression
naturally approaches, but need not exactly reach, unity in the absence of
repressor, capturing biological realities such as polymerase saturation and
stochastic promoter switching that the Hill function's fixed maximum cannot
represent.

The systematic parameter estimation strategy developed here bridges existing
Hill-linear hybrid formulations and pure logistic reformulations by
simultaneously matching basal expression rates and local activation slopes. The
biologically motivated threshold criterion $\theta_B = g_A$---centering
the sigmoid's inflection point precisely where the weighted cross-activation
signal $g_{AB}B$ equals basal production---yields the explicit closed-form
relationships $\kappa_1 = 4g_A$,
$\theta_B = g_A$, and $\lambda_1 = \ln 3/g_A$ for the weighted logistic
activation, alongside $\lambda = n/\theta$ for logistic self-repression. The
logistic steepness parameter $\lambda$ further provides a continuous tuning
mechanism between the linear regime ($\lambda \ll 1$) and the sharp switch-like
regime ($\lambda \gg 1$), enabling smooth interpolation across dynamical
behaviors ranging from genetic oscillations and bistable switches to graded
developmental patterning.

Future work will extend the approach to direct numerical simulation of complex
dynamics and extreme events in the logistic reformulations, to multi-gene
regulatory networks, to stochastic
formulations accounting for intrinsic noise in low-copy-number regimes, and
to experimental validation in optogenetic platforms. By replacing Hill functions
with their logistic counterparts while preserving sigmoidal dynamics, researchers
can build on decades of accumulated Hill-based modeling intuition while gaining
the analytical tractability, numerical reliability, and biological fidelity that
demanding applications in synthetic biology, metabolic engineering, and
therapeutic cell design require.

\section{Declarations}

\subsection*{Availability of Data, Materials and Code}
All derivations and numerical results are reproducible from the equations and
parameters given in the manuscript. The R scripts that produced every
numerical value reported in the paper are provided as supplementary material:
\begin{itemize}
\item \texttt{equilibrium\_linear\_additive.R} --- Newton--Raphson computation of
the delay-free equilibrium $E^* = (\bar A,\bar B)$ for the linear additive
formulation~\eqref{eq:A_logistic}--\eqref{eq:B_logistic}, including the
logistic function values $f_A^-,\,f_B^-$ at $E^*$.
\item \texttt{equilibrium\_weighted\_logistic.R} --- analogous Newton--Raphson
solver for the weighted logistic
formulation~\eqref{eq:A_logistic_basal}--\eqref{eq:B_logistic_basal},
returning $E^*_{\mathrm{wl}}$ and the four logistic factors
$f_1^+,\,f_2^+,\,f_3^-,\,f_4^-$.
\item \texttt{hopf\_bifurcation\_solver.R} --- joint solution of the
two-equation transcendental Hopf system~\eqref{eq:hopf_system} via grid-seeded
Newton iterations, returning $(\omega_c,\,\tau_c)$, the transversality
coefficient $d\,\mathrm{Re}(\mu)/d\tau$, and the higher-order Hopf points
$\tau_c^{(k)}$ for both formulations.
\end{itemize}
The scripts use only the base R distribution together with the
\texttt{nleqslv} package (CRAN). Each script is self-contained and prints
all numerical values reported in Sections~\ref{sec:equilibrium}
and~\ref{sec:stability} to the console; running the three scripts in any
order reproduces every quantitative result of the paper.

\subsection*{Competing Interests}
The author declares no competing interests.

\subsection*{Funding}
Not applicable.

\subsection*{Ethics Approval and Consent to Participate}
Not applicable.

\subsection*{Authors' Contributions}
Single-author study: all aspects conceived, derived, computed, and written by the author.

\subsection*{Acknowledgements}
Not applicable.

\bibliographystyle{elsarticle-num}
\bibliography{mybibfile}

\begin{thebibliography}{10}
\expandafter\ifx\csname url\endcsname\relax
  \def\url#1{\texttt{#1}}\fi
\expandafter\ifx\csname urlprefix\endcsname\relax\def\urlprefix{URL }\fi
\expandafter\ifx\csname href\endcsname\relax
  \def\href#1#2{#2} \def\path#1{#1}\fi

\bibitem{farcot2019chaos}
E.~Farcot, S.~Best, R.~Edwards, I.~Belgacem, X.~Xu, P.~Gill, Chaos in a ring
  circuit, Chaos: An Interdisciplinary Journal of Nonlinear Science 29~(4)
  (2019) 043103.

\bibitem{belgacem2025glass}
I.~Belgacem, R.~Edwards, {\'E}.~Farcot, Computer-aided analysis of
  high-dimensional {Glass} networks: periodicity, chaos, and bifurcations in a
  ring circuit, Chaos: An Interdisciplinary Journal of Nonlinear
  ScienceAccepted 29 January 2025 (may 2025).
\newblock \href {https://doi.org/10.1063/5.0243955}
  {\path{doi:10.1063/5.0243955}}.

\bibitem{belgacem2018bmb}
I.~Belgacem, S.~Casagranda, E.~Grac, D.~Ropers, J.-L. Gouz{\'e}, Reduction and
  stability analysis of a transcription-translation model of {RNA} polymerase,
  Bulletin of Mathematical Biology 80 (2018) 294--318.
\newblock \href {https://doi.org/10.1007/s11538-017-0372-4}
  {\path{doi:10.1007/s11538-017-0372-4}}.

\bibitem{belgacem2014cdcRNA}
I.~Belgacem, E.~Grac, D.~Ropers, J.-L. Gouz{\'e}, Stability analysis of a
  reduced transcription-translation model of {RNA} polymerase, in: Proceedings
  of the {IEEE} 53rd Conference on Decision and Control ({CDC}), IEEE, Los
  Angeles, California, USA, 2014, pp. 3924--3929.
\newblock \href {https://doi.org/10.1109/CDC.2014.7039999}
  {\path{doi:10.1109/CDC.2014.7039999}}.

\bibitem{belgacem2014med}
I.~Belgacem, J.-L. Gouz{\'e}, Mathematical study of the global dynamics of a
  concave gene expression model, in: Proceedings of the 22nd Mediterranean
  Conference on Control and Automation ({MED}), IEEE, Palermo, Italy, 2014, pp.
  1341--1346.
\newblock \href {https://doi.org/10.1109/MED.2014.6961562}
  {\path{doi:10.1109/MED.2014.6961562}}.

\bibitem{belgacem2013cdc}
I.~Belgacem, J.-L. Gouz{\'e}, Stability analysis and reduction of gene
  transcription models, in: Proceedings of the {IEEE} 52nd Annual Conference on
  Decision and Control ({CDC}), IEEE, Florence, Italy, 2013, pp. 2691--2696.
\newblock \href {https://doi.org/10.1109/CDC.2013.6760289}
  {\path{doi:10.1109/CDC.2013.6760289}}.

\bibitem{belgacem2013cab}
I.~Belgacem, J.-L. Gouz{\'e}, Analysis and reduction of transcription
  translation coupled models for gene expression, in: Proceedings of the 12th
  {IFAC} Symposium on Computer Applications in Biotechnology ({CAB}), Elsevier,
  Mumbai, India, 2013, pp. 42--47.
\newblock \href {https://doi.org/10.3182/20131216-3-IN-2044.00012}
  {\path{doi:10.3182/20131216-3-IN-2044.00012}}.

\bibitem{belgacem2013acta}
I.~Belgacem, J.-L. Gouz{\'e}, Global stability of enzymatic chain of full
  reversible {Michaelis--Menten} reactions, Acta Biotheoretica 61~(3) (2013)
  425--436.
\newblock \href {https://doi.org/10.1007/s10441-013-9195-3}
  {\path{doi:10.1007/s10441-013-9195-3}}.

\bibitem{belgacem2012ifac}
I.~Belgacem, J.-L. Gouz{\'e}, Global stability of full open reversible
  {Michaelis--Menten} reactions, in: Proceedings of the 8th {IFAC} Symposium on
  Advanced Control of Chemical Processes, Elsevier, Singapore, 2012, pp.
  591--596.
\newblock \href {https://doi.org/10.3182/20120710-4-SG-2026.00039}
  {\path{doi:10.3182/20120710-4-SG-2026.00039}}.

\bibitem{belgacem2019ecc}
I.~Belgacem, H.~Bensalah, B.~Cherki, R.~Edwards, The probabilistic convolution
  regularization of {Zeno} hybrid systems, in: Proceedings of the 18th European
  Control Conference ({ECC}), IEEE, Naples, Italy, 2019.
\newblock \href {https://doi.org/10.23919/ECC.2019.8795626}
  {\path{doi:10.23919/ECC.2019.8795626}}.

\bibitem{belgacem2020cdc}
I.~Belgacem, J.-L. Gouz{\'e}, R.~Edwards, Control of negative feedback loops in
  genetic networks, in: Proceedings of the 59th {IEEE} Conference on Decision
  and Control ({CDC}), IEEE, Jeju Island, Republic of Korea, 2020.
\newblock \href {https://doi.org/10.1109/CDC42340.2020.9304088}
  {\path{doi:10.1109/CDC42340.2020.9304088}}.

\bibitem{chambon2020automatica}
L.~Chambon, I.~Belgacem, J.-L. Gouz{\'e}, Qualitative control of undesired
  oscillations in a genetic negative feedback loop with uncertain measurements,
  Automatica 112 (2020) 108642.
\newblock \href {https://doi.org/10.1016/j.automatica.2019.108642}
  {\path{doi:10.1016/j.automatica.2019.108642}}.

\bibitem{ismail2025logistic}
B.~Ismail, \href{https://arxiv.org/abs/2512.14325}{Exploring logistic functions
  as robust alternatives to hill functions in genetic network modeling},
  submitted December 16, 2025 (2025).
\newblock \href {http://arxiv.org/abs/2512.14325} {\path{arXiv:2512.14325}}.
\newline\urlprefix\url{https://arxiv.org/abs/2512.14325}

\bibitem{vinoth2025extreme}
S.~Vinoth, S.~L. Kingston, S.~Srinivasan, S.~Kumarasamy, T.~Kapitaniak, Extreme
  events in gene regulatory networks with time-delays, Scientific Reports
  15~(1) (2025) 13064.

\bibitem{madar2011negative}
D.~Madar, E.~Dekel, A.~Bren, U.~Alon, Negative auto-regulation increases the
  input dynamic-range of the arabinose system of escherichia coli, BMC systems
  biology 5~(1) (2011) 111.

\bibitem{oehler1994quality}
S.~Oehler, M.~Amouyal, P.~Kolkhof, B.~von Wilcken-Bergmann, B.~M{\"u}ller-Hill,
  Quality and position of the three lac operators of e. coli define efficiency
  of repression., The EMBO journal 13~(14) (1994) 3348--3355.

\bibitem{jacob1961genetic}
F.~Jacob, J.~Monod, Genetic regulatory mechanisms in the synthesis of proteins,
  Journal of molecular biology 3~(3) (1961) 318--356.

\end{thebibliography}

\end{document}